\ifdefined\submit
\OneAndAHalfSpacedXI % current default line spacing
\usepackage{graphicx}
\usepackage{multirow}
\usepackage{hhline}
\usepackage[shortlabels]{enumitem}
\usepackage[small, margin=1cm]{caption}

\usepackage{appendix}
\usepackage{color}
\definecolor{strcolor}{rgb}{0.6, 0.2, 0.6}
\definecolor{commentcolor}{rgb}{0.3125, 0.5, 0.3125}
\definecolor{keycol}{rgb}{0, 0, 1}
\usepackage{xcolor}
 \newcommand{\alert}[1]{\textcolor{magenta}{#1}}
\usepackage{bbm}
\usepackage{mathrsfs}
\usepackage{mathtools}
\usepackage[breaklinks=true,pdfstartview=FitH]{hyperref}
\usepackage[capitalize,nameinlink]{cleveref}
\crefname{assumption}{Assumption}{Assumptions}

\usepackage{listings}
\usepackage{url}

\newcommand {\bea}{\begin{eqnarray}}
	\newcommand {\eea}{\end{eqnarray}}

\newcommand{\vmuquad}{v_{\mu}^{\rm quad}}
\newcommand{\smuquad}{S_{\mu}}
\newcommand{\vmuentrop}{v_{\mu}^{\rm ent}}

\renewcommand{\Re}{\mathbb{R}}

\newcommand{\ds}{\displaystyle}

\newcommand{\vmu}{v_{\mu}}
\newcommand{\gmu}{g_{\mu}}
\newcommand{\phimu}{\phi_{\mu}}
\newcommand{\gtildemu}{\tilde{g}_{\mu}}
\newcommand{\Ymu}{Y_{\mu}}
\newcommand{\citen}[2]{\textnormal{\textsc{\cite[#1]{#2}}}}

\newcommand{\jhsolved}[1]{}

\def \rla{\rangle}
\def \lla{\langle}

\DeclareMathOperator*{\diam}{diam}

\DeclareMathOperator*{\cl}{cl}
\DeclareMathOperator*{\dist}{dist}
\DeclareMathOperator*{\interior}{int}
\DeclareMathOperator*{\co}{co}
\DeclareMathOperator*{\vol}{vol}

\DeclareMathOperator*{\dom}{dom}
\newcommand{\norm}[1]{{\left\|{#1}\right\|}}

\def\blot{\quad \mbox{$\vcenter{ \vbox{ \hrule height.4pt
				\hbox{\vrule width.4pt height.9ex \kern.9ex \vrule width.4pt}
				\hrule height.4pt}}$}}

\usepackage{natbib}
\bibpunct[, ]{(}{)}{,}{a}{}{,}%
\TheoremsNumberedThrough     % Preferred (Theorem 1, Lemma 1, Theorem 2)
\ECRepeatTheorems

\EquationsNumberedThrough    % Default: (1), (2), ...
\gdef\AQ#1{}
\gdef\CQ#1{}
\else

\documentclass[letter,11pt]{article}
\usepackage{amsthm}
\usepackage{amssymb,amscd}
\usepackage{wrapfig}
\usepackage[pdftex]{graphicx}
\usepackage[latin1]{inputenc}
\usepackage[framemethod=tikz]{mdframed}
\usepackage[ruled,vlined]{algorithm2e}
\usepackage{epsfig,epstopdf}
\usepackage{mathrsfs}
\usepackage{mathtools}
\usepackage{bbm}
\usepackage{seqsplit}
\usepackage{multirow}
\usepackage{caption}   
\usepackage{float} 
\usepackage{booktabs}
\usepackage{bigstrut}
\usepackage{authblk}
\usepackage{tcolorbox}

\usepackage{natbib}
\usepackage[breaklinks=true,pdfstartview=FitH,pagebackref=true]{hyperref}
\usepackage{footnotebackref}
\usepackage{tablefootnote}
\usepackage[shortlabels]{enumitem}
\usepackage{color}
\usepackage{xcolor}
\usepackage{subcaption}

\newcommand{\ds}{\displaystyle}

\newcommand{\vmu}{v_{\mu}}
\newcommand{\gmu}{g_{\mu}}
\newcommand{\phimu}{\phi_{\mu}}
\newcommand{\gtildemu}{\tilde{g}_{\mu}}
\newcommand{\Ymu}{Y_{\mu}}
\newcommand{\vmuquad}{v_{\mu}^{\rm quad }}
\newcommand{\smuquad}{S_{\mu}}
\newcommand{\vmuentrop}{v_{\mu}^{\rm ent}}
\newcommand{\citen}[2]{\textnormal{\cite[#1]{#2}}}
 \newcommand{\alert}[1]{#1}

\newcommand{\jhsolved}[1]{}

\def \rla{\rangle}
\def \lla{\langle}

\renewcommand{\Re}{{\rm I}\! {\rm R}}

\DeclareMathOperator*{\diam}{diam}

\DeclareMathOperator*{\argmin}{arg\,min}

\DeclareMathOperator*{\cl}{cl}
\DeclareMathOperator*{\dist}{dist}
\DeclareMathOperator*{\interior}{int}
\DeclareMathOperator*{\co}{co}
\DeclareMathOperator*{\vol}{vol}

\DeclareMathOperator*{\dom}{dom}
\newcommand{\norm}[1]{{\left\|{#1}\right\|}}

\usepackage[margin=1in]{geometry}
\usepackage[capitalize,nameinlink]{cleveref}
\crefname{assumption}{Assumption}{Assumptions}
\crefformat{equation}{\textup{#2(#1)#3}}
\crefrangeformat{equation}{\textup{#3(#1)#4--#5(#2)#6}}
\crefmultiformat{equation}{\textup{#2(#1)#3}}{ and \textup{#2(#1)#3}}
{, \textup{#2(#1)#3}}{, and \textup{#2(#1)#3}}
\crefrangemultiformat{equation}{\textup{#3(#1)#4--#5(#2)#6}}%
{ and \textup{#3(#1)#4--#5(#2)#6}}{, \textup{#3(#1)#4--#5(#2)#6}}{, and \textup{#3(#1)#4--#5(#2)#6}}

\Crefformat{equation}{#2Equation~\textup{(#1)}#3}
\Crefrangeformat{equation}{Equations~\textup{#3(#1)#4--#5(#2)#6}}
\Crefmultiformat{equation}{Equations~\textup{#2(#1)#3}}{ and \textup{#2(#1)#3}}
{, \textup{#2(#1)#3}}{, and \textup{#2(#1)#3}}
\Crefrangemultiformat{equation}{Equations~\textup{#3(#1)#4--#5(#2)#6}}%
{ and \textup{#3(#1)#4--#5(#2)#6}}{, \textup{#3(#1)#4--#5(#2)#6}}{, and \textup{#3(#1)#4--#5(#2)#6}}

\crefdefaultlabelformat{#2\textup{#1}#3}

\numberwithin{equation}{section}

\newtheorem{theorem}{Theorem}[section]
\newtheorem{definition}[theorem]{Definition}

\newtheorem{proposition}[theorem]{Proposition}
\newtheorem{lemma}[theorem]{Lemma}
\newtheorem{assumption}[theorem]{Assumption}
\theoremstyle{definition} 	
\newtheorem{example}[theorem]{Example}
\newtheorem{remark}[theorem]{Remark}

\def\TheTitle{\alert{Theoretical smoothing frameworks for nonsmooth simple bilevel problems}}

\ifpdf
\hypersetup{
	pdftitle={\TheTitle},
	pdfauthor={Jan Harold Alcantara and Akiko Takeda}}
\fi

\title{\TheTitle}
\date{\today}
\author{Jan Harold Alcantara\thanks{\url{janharold.alcantara@riken.jp}.  Center for Advanced Intelligence Project, RIKEN, Tokyo, Japan.}
\qquad Akiko Takeda\thanks{\url{takeda@mist.i.u-tokyo.ac.jp}.
	Department of Mathematical Informatics, Graduate School of Information Science and Technology, University of Tokyo,
	Tokyo, Japan Center for Advanced Intelligence Project, RIKEN, Tokyo, Japan}
}
\fi 
\begin{document}
	
\ifdefined\submit 
	%%%%%%%%%%%%%%%%
	
%	\AIA
% \setcounter{page}{1} %
% \VOL{00}%
% \NO{0}%
% \MONTH{Xxxxx}%
% \YEAR{2017}%
% \FIRSTPAGE{1}%
% \LASTPAGE{16}%
% \FIRSTPAGEAIA{1}%
% \LASTPAGEAIA{16}%
\def\COPYRIGHTHOLDER{INFORMS}%
\def\COPYRIGHTYEAR{2017}%
\def\DOI{\fontsize{7.5}{9.5}\selectfont\sf\bfseries\noindent https://doi.org/10.1287/opre.2017.1714\CQ{Word count = 9740}}
%\def\RECEIVED{November 1, 2016}
%\def\REVISED{June 22, 2017; October 6, 2017}
%\def\ACCEPTED{November 15, 2017}
% \PUBONLINEAIA{}

	\RUNAUTHOR{Alcantara and Takeda} %

	\RUNTITLE{Theoretical smoothing frameworks for nonsmooth simple bilevel problems}

\TITLE{\alert{Theoretical smoothing frameworks for nonsmooth simple bilevel problems}}


	\ARTICLEAUTHORS{
%		\AUTHOR{Jianzhe Zhen,\textsuperscript{a,*} Dick den
%		Hertog,\textsuperscript{a} Melvyn Sim\textsuperscript{b}} 
%\AFF{$^{a}$Department of Econometrics and Operations Research,
%Tilburg University; $^{b}$NUS Business School, National University of
%Singapore}

\AUTHOR{Jan Harold Alcantara}
\AFF{ Center for Advanced Intelligence Project, RIKEN, Tokyo, Japan}

\AUTHOR{Akiko Takeda}
\AFF{	Department of Mathematical Informatics, Graduate School of Information Science and Technology, University of Tokyo,
	Tokyo, Japan and \\ Center for Advanced Intelligence Project, RIKEN, Tokyo, Japan}

%\AUEXTRA{$^{*}$Corresponding author}

%\AFFmail{{\bf Contact:} j.zhen@tilburguniversity.edu,
%d.denhertog@tilburguniversity.edu,\\			melvynsim@nus.edu.sg}%
}
	 % end of the block
	
%\ARTICLEINFO{\textbf{Received:} November 1, 2016\\ \textbf{Revised:} June 22, 2017; October 6, 2017\\ \textbf{Accepted:} November 15, 2017\\ \textbf{Published Online in Articles in Advance:}}

	\ABSTRACT{Bilevel programming has recently received a great deal of attention due to its abundant applications in many areas. The optimal value function approach provides a useful reformulation of the bilevel problem, but its utility is often limited due to the nonsmoothness of the value function even in cases when the associated lower-level function is smooth. In this paper, we present two smoothing strategies for the value function associated with lower-level functions that are not necessarily smooth but are Lipschitz continuous. The first method employs quadratic regularization for partially convex lower-level functions, while the second utilizes entropic regularization for general lower-level objective functions. Meanwhile, the property known as gradient consistency is crucial in ensuring that a designed smoothing algorithm is globally subsequentially convergent to stationary points of the value function reformulation. With this motivation, we prove that the proposed smooth approximations satisfy the gradient consistent property under certain conditions on the lower-level function.}

%\FUNDING{The research of the first author is supported by NWO Grant 613.001.208. The third author acknowledges the funding support from the Singapore Ministry of Education Social Science Research Thematic Grant MOE2016-SSRTG-059.}

%\SUBJECTCLASS{\AQ{Please confirm subject classifications.}Fourier-Motzkin elimination; adjustable robust optimization; linear decision rules; redundant constraint identification.}

\AREAOFREVIEW{Continuous Optimization.}

%{\CQ{Kindly provide the keywords.}}

	%%%%%%%%%%%%%%%%%%%%%%%%%%%%%%%%%%%%%%%%%%%%%%%%%%%%%%%%%%%%%%%%%%%%%%
	
	% Samples of sectioning (and labeling) in OPRE
	% NOTE: (1) \section and \subsection do NOT end with a period
	%       (2) \subsubsection and lower need end punctuation
	%       (3) capitalization is as shown (title style).
	%
	%\section{Introduction.}\label{intro} %%1.
	%\subsection{Duality and the Classical EOQ Problem.}\label{class-EOQ} %% 1.1.
	%\subsection{Outline.}\label{outline1} %% 1.2.
	%\subsubsection{Cyclic Schedules for the General Deterministic SMDP.}
	%  \label{cyclic-schedules} %% 1.2.1
	%\section{Problem Description.}\label{problemdescription} %% 2.
	% Text of your paper here

\maketitle
\else

\maketitle 
\begin{abstract}
	Bilevel programming has recently received a great deal of attention due to its abundant applications in many areas. The optimal value function approach provides a useful reformulation of the bilevel problem, but its utility is often limited due to the nonsmoothness of the value function even in cases when the associated lower-level function is smooth. In this paper, we present two smoothing strategies for the value function associated with lower-level functions that are not necessarily smooth but are Lipschitz continuous. The first method employs quadratic regularization for partially convex lower-level functions, while the second utilizes entropic regularization for general lower-level objective functions. Meanwhile, the property known as gradient consistency is crucial in ensuring that a designed smoothing algorithm is globally subsequentially convergent to stationary points of the value function reformulation. With this motivation, we prove that the proposed smooth approximations satisfy the gradient consistent property under certain conditions on the lower-level function.\\
	
	%   	\ifdefined\submit
	%   	\keywords{fixed point algorithm; proximal methods; 
	%   		alternating projections; averaged projections; linear complementarity 
	%   		problem; union convex set; nonconvex feasibility problems; nonconvex optimization; global convergence }
	%   	
	%   	\else 
\end{abstract}
\fi

\section{Introduction}
Consider the bilevel optimization problem
\begin{align}
\label{eq:bilevelproblem}
\tag{BP}
\begin{array}{rl}
\ds \min _{(x,y)\in X\times Y}~ & f(x,y) \\
\mbox{s.t.} ~& \ds y \in \argmin_{\bar{y}\in Y} g(x,\bar{y})
\end{array}
\end{align}
where $f,g:\Re^n\times \Re^m \to \Re$ are Lipschitz continuous but possibly nonsmooth functions, and $X,Y$ are closed and convex subsets of $\Re^n$ and $\Re^m$, respectively.\jhsolved{Consider as well the following: Restrict $f$ and $g$ to the open set $O_X\times O_Y$ that contains $X\times Y$. $f$ has to be continuous on closure of $O_X\times O_Y$ (dont forget notation for closure) and $C^1$ on $O_X\times O_Y$ and $g$ is Lipschitz continuous on $O_X\times O_Y$ (need open sets)... Convexity assumptions?? Where is $V$ and $S$ defined? Also, $f$ does not need to be differentiable....} The minimization of $g(x,\cdot)$ over $Y$ is known as the \emph{lower-level problem}, while the minimization of $f$ is called the \emph{upper-level problem}. The bilevel optimization problem \eqref{eq:bilevelproblem} has recently gained significant attention due to its numerous applications, including hyperparameter optimization \citep{franceschi18a,OTKW21}, meta-learning \citep{franceschi18a}, and Stackelberg games \citep{von2010market}, among others. We refer the reader to \citep{dempe2015bilevel,dempe2020bilevel} and the references therein for more applications of bilevel optimization.

Traditional methods for solving the bilevel problem \eqref{eq:bilevelproblem} involve reformulation or relaxation of \eqref{eq:bilevelproblem} as a one-level problem, thereby permitting the use of solution methods from nonlinear programming \citep{dempe2015bilevel}. A widely used approach is to replace the lower-level minimization problem by its optimality conditions. This yields an equivalent problem if $g(x,\cdot)$ is convex for each $x\in X$. Among the recent trends that gained considerable attention is the case when $f$ is continuously differentiable, $g$ is a twice continuously differentiable, $g(x,\cdot)$ is strongly convex for all $x\in X$, and $Y=\Re^m$. In this case, the lower-level problem of \eqref{eq:bilevelproblem} defines a smooth map $x\mapsto y(x)$ that corresponds to the unique minimizer of $g(x,\cdot)$. In turn, \eqref{eq:bilevelproblem} simplifies to the minimization of a smooth function $f(x,y(x))$ over $x\in X$, which can be handled by typical smooth optimization methods such as the projected gradient descent algorithm. Indeed, this framework has been the basis of some recent works in the literature \citep{GhadimiWang18,ji21c,pedregosa16}, but the main drawback is that this approach is heavily dependent on the smoothness of $f$ and $g$, as well as the strong convexity of $g(x,\cdot)$. In many applications, such smoothness requirements are not met \citep{OTKW21}. Moreover, we note that when $g(x,\cdot)$ is nonconvex, a local optimal solution of \eqref{eq:bilevelproblem} may not even correspond to a stationary point of its reformulation via optimality conditions of the lower-level problem \citep{Mirrlees99}.

An alternative stream of research, which is the direction pursued in this paper, employs the {value function reformulation} approach. The \textit{value function} associated with the lower-level problem of \eqref{eq:bilevelproblem} is defined as
\begin{equation}
v(x) \coloneqq \min _{y\in Y}~ g(x,y). \label{eq:valuefunction}
\end{equation}

With this, the bilevel problem \eqref{eq:bilevelproblem} is equivalent to the nonlinear programming problem
\begin{align}
\begin{array}{rl}
\ds \min _{(x,y)\in X\times Y}~ & f(x,y) \\
\text{s.t.} ~& \ds g(x,y) - v(x) \leq 0,
\end{array}
\tag{VFP}
\label{eq:bilevel_nlp_formulation}
\end{align}
which was first proposed by \citet{Outrata90}. 
Note that the reformulation \eqref{eq:bilevel_nlp_formulation} is independent of any smoothness and convexity assumptions on $g$, so that this framework is theoretically applicable to any bilevel optimization problem. However, one difficulty that comes with this flexible reformulation is determining under what conditions local solutions of \eqref{eq:bilevel_nlp_formulation} correspond to its stationary points. In nonlinear programming, these conditions pertain to constraint qualifications (CQs), and unfortunately, usual CQs such as the Mangasarian-Fromovitz constraint qualification (MFCQ) are violated by feasible points of \eqref{eq:bilevel_nlp_formulation} \citep{YeZhu95}. One promising alternative is to consider an approximate bilevel program 
\begin{align}
\begin{array}{rl}
\ds \min _{(x,y)\in X\times Y}~ & f(x,y) \\
\text{s.t.} ~& \ds g(x,y) - v(x) \leq \epsilon.
\end{array}
\tag{$\text{VFP}_{\epsilon}$}
\label{eq:bilevel_nlp_formulation_epsilon}
\end{align}
where $\epsilon > 0$ \cite{LinXuYe14,lu2023slm,YeYuanZengZhang23}. The above perturbed problem is a reasonable approximation of \eqref{eq:bilevel_nlp_formulation} as its local and global solutions are arbitrarily close to the solution set of \eqref{eq:bilevel_nlp_formulation} \citep{LinXuYe14,YeYuanZengZhang23}.
%; see \cite[Proposition 6]{YeYuanZengZhang23}. 
Moreoever, unlike \eqref{eq:bilevel_nlp_formulation}, (nonsmooth) MFCQ is automatically satisfied by feasible points of \eqref{eq:bilevel_nlp_formulation_epsilon}, which is a sufficient constraint qualification for local solutions of \eqref{eq:bilevel_nlp_formulation_epsilon} to be stationary \citep{LamparielloSagratella20,LinXuYe14,YeYuanZengZhang23}. 
%. Using the above approximation, the (nonsmooth) MFCQ is more likely to hold (see discussions in \cite{LinXuYe14}). In fact, if $g$ is smooth, $g(x,\cdot)$ is a convex function for any $x\in X$ and $Y$ is a convex set, then the (smooth) MFCQ always holds at any feasible point of \eqref{eq:bilevel_nlp_formulation_epsilon}. See also \cite[Definition 2]{LamparielloSagratella20}.

While the approximate bilevel problem \eqref{eq:bilevel_nlp_formulation_epsilon} addresses the CQ issue, another hurdle is the need for an algorithmic framework that is capable of obtaining stationary points. Solution methods for solving nonlinear programming problems often hinge on the assumption that both the objective and constraint functions involved are smooth. Unfortunately, the value function reformulation \eqref{eq:bilevel_nlp_formulation} of the bilevel problem introduces a challenge, since the value function $v$ defined by \eqref{eq:valuefunction}, in general, is nonsmooth even when $g$ is a smooth function (see \cref{thm:danskin}). In the absence of smoothness, a popular adopted strategy involves employing smooth approximations, which have been extensively utilized in many nonsmooth optimization problems, complementarity problems and variational inequalities; see \cite{ChenXiaojun12} and the references therein. As a consequence, it is natural to expect that smooth approximations for the functions $f$ and $g$ in the bilevel problem \eqref{eq:bilevelproblem} may be readily available, since these often correspond to loss functions with a nonsmooth regularizer in many applications; for instance, see \cite{ChenXiaojun12}. However, as previously mentioned, the smoothness of $g$ is in general not adequate to ensure the smoothness of $v$. Consequently, employing traditional methods\jhsolved{Fixed this part.} for the nonlinear programming reformulation \eqref{eq:bilevel_nlp_formulation} or \eqref{eq:bilevel_nlp_formulation_epsilon} necessitates smoothing techniques for the value function. Unfortunately, to date, smooth approximations of the value functions have not been thoroughly investigated yet in the literature. One work that pursued this direction in the context of bilevel optimization is studied by \cite{LinXuYe14}. The authors adopted the entropic regularization approach originally formulated in \citep{FangWu96}, which yields an approximation (in the sense of \cref{defn:smoothing}) that is valid for a smooth function $g$ and a compact set $Y$.  \textit{Gradient consistency} (see \cref{defn:gradconsistency}) of the proposed smoothing family is also proved, which is an essential property to guarantee convergence to stationary points. Apart from this work, there has been no other existing works on smoothing methods for value function reformulation of bilevel problems, to the best of our knowledge.

\alert{\subsubsection*{Contributions of this work} In this paper, we present a theoretical framework for constructing and analyzing \emph{smooth approximations} $v_\mu$ of the {value function} \eqref{eq:valuefunction}, where $g$ may be a nonsmooth function. Our primary goal is to design smoothing families $\{v_\mu\}$ that not only approximate $v$ well but also satisfy the crucial property of {gradient consistency}, thereby establishing a foundation for convergence analysis of smoothing algorithms based on the value function approach to bilevel optimization.}

\alert{We focus on the {value function reformulation} of the bilevel problem and propose two smoothing techniques applicable to Lipschitz continuous functions $g$. The first is a quadratic regularization method, suitable when $g$ is convex in $y$. The second is a more general entropic regularization approach, which does not require convexity and extends the framework of \citep{FangWu96,LinXuYe14} to derive smooth approximations of the value function over compact sets $Y$ when $g$ is nonsmooth. We also discuss an extension to unbounded closed sets.}

\alert{As mentioned above, gradient consistency plays a central role in guaranteeing convergence of smoothing algorithms to stationary points (see \cref{prop:stationary_smoothing}). Accordingly, we provide sufficient conditions under which the proposed smooth approximations possess this property. A key technical contribution is the development of {Danskin-type theorems} for value functions over arbitrary closed sets $Y$, which serve as the foundation for our analysis.}

\alert{\subsubsection*{Organization of the paper} This paper is organized as follows. In \cref{sec:prelim}, we review essential concepts from nonsmooth and variational analysis, gradient consistency, and optimality conditions for the value function reformulation. In \cref{sec:quadratic}, we present our first main result: a smoothing approach based on quadratic regularization, along with its gradient consistency properties. The second technique, entropic regularization, is developed in \cref{sec:entropy}, where we also analyze its theoretical properties. Finally, concluding remarks and directions for future work are provided in \cref{sec:conclusion}.}

\section{Mathematical Preliminaries}\label{sec:prelim}
Throughout the paper, $B_{\delta}(x)$ denotes an open ball centered at $x$ with radius $\delta>0$. The \textit{interior}, \textit{closure} and \textit{convex hull} of a set $X$ are denoted by $\interior (X)$, $\cl (X)$ and $\co (X)$, respectively. The \textit{projection} onto $X$ and the \textit{distance function} on $X$ are given respectively by 
\[ \Pi_X(z) \coloneqq \argmin \{\|x-z\| : x\in X \},\]
and 
\[ \dist (z,X) \coloneqq \min \{ \|x-z\| : x\in X\}, \]
\alert{where $\norm{\cdot}$ denotes the Euclidean norm on $\Re^n$. The indicator function of the set $X$, denoted by $\delta_X$, is defined as $\delta_X(x) = 0$ if $x \in X$, and $\delta_X(x) = \infty$ otherwise.}
We will always use $O_X$ and $O_Y$ to denote an open set containing the closed sets $X$ and $Y$, respectively.  \alert{The Lebesgue measure of $X$ is denoted by $\vol (X)$.} For a set $S\subseteq \Re^n\times \Re^m$, we define 
\begin{equation}
\pi_x (S) \coloneqq \{x' \in \Re^n : \exists y' \in \Re^m ~\text{s.t. } (x',y')\in S \}.
\label{eq:pi_x}
\end{equation}

\subsection{Nonsmooth and variational analysis}\label{sec:prelim_nonsmooth}

We recall some important concepts and results from nonsmooth and variational analysis. Most of these materials can be found in \cite[Chapter 2]{Clarke83} and \cite{RW98}.

Let $\phi$ be a Lipschitz continuous function on an open set $O\subseteq \Re^n$. The \emph{Clarke generalized directional derivative} of $\phi$ at $\bar{x}\in O$ in the direction $d$, denoted by $\phi^{\circ}(x;d)$, is defined as
\[\phi^{\circ}(\bar{x};d) = \limsup_{x\to \bar{x}, t\searrow 0} \frac{\phi (x+td) - \phi(x)}{t} ,\]
and the (usual) \emph{directional derivative} of $\phi$ at $\bar{x}$ in the direction $d$ is 
\[\phi'(\bar{x};d) = \lim_{t\searrow 0} \frac{\phi (\bar{x}+td) - \phi(\bar{x})}{t} . \]
We say that $\phi$ is  at $\bar{x}$ \emph{Clarke regular} if for all $d\in \Re^n$, $\phi'(\bar{x};d)$ exists and $\phi'(\bar{x};d) = \phi^{\circ}(\bar{x};d)$. For instance, if $\phi$ is differentiable at $\bar{x}$ or if $\phi$ is convex, then it is Clarke regular at $\bar{x}$. 

The \emph{Clarke generalized gradient} of $\phi$ at $\bar{x}$, denoted by $\partial \phi(\bar{x})$, is given by 
\[\partial \phi (\bar{x}) \coloneqq \{ \xi\in \Re^n : \phi^{\circ}(\bar{x};d)\geq \lla \xi, d\rla ~\forall d\in\Re^n  \} ,\]
which is a compact and convex set. When $\phi$ is convex, this coincides with the subdifferential in the sense of convex analysis:
\[ \partial^c \phi(\bar{x}) \coloneqq \{ \xi \in \Re^n: \phi(x) \geq \phi(\bar{x}) + \lla \xi, x-\bar{x}\rla ~\forall x\in\Re^n \} .\]
The Clarke generalized directional derivative $\phi^{\circ}(\bar{x};\cdot)$ can be expressed as the support function of $\partial\phi(\bar{x})$, that is, 
\begin{equation}
\phi^{\circ}(\bar{x};d) = \sigma_{\partial \phi (\bar{x})} ( d) \coloneqq \max \{ \lla d, \xi \rla : \xi \in \partial \phi(\bar{x})\},
\label{eq:maxformula}
\end{equation}
and so when $\phi$ is convex, $\phi'(\bar{x};d) = \sigma_{\partial^c \phi(\bar{x})}(d)$. 

As a set-valued mapping, $\partial \phi (\cdot)$ is upper semicontinuous at any $\bar{x}\in O$. That is, given any $\varepsilon>0$, there exists $\delta>0$ such that\jhsolved{$B_0(1)$ is now replaced with $B_1(0)$.} 
\[\partial \phi(x)\subseteq \partial \phi (\bar{x})+ \varepsilon B_1(0) \quad \forall x\in B_{\delta}(\bar{x}). \]
By a standard argument using Heine-Borel theorem and noting the boundedness of $\partial \phi (x)$ for any $x\in O$, note that $\bigcup_{x\in \Omega} \partial \phi (x)$ is a bounded set for any compact subset $\Omega$ of $O$.

For the Clarke subdifferential of the value function, we have the following celebrated result due to Danskin (see \cite{Danskin67} and \cite[Ex. 9.13]{Clarke98}).

\begin{theorem}[Danskin's Theorem]
	Let $g:\Re^n \times Y\to \Re$ be given, where $Y$ is a compact subset of $\Re^m$. Suppose that for a neighborhood $\Omega$ of $\bar{x}\in \Re^n$, the derivative $\nabla_x g(x,y)$ exists and is continuous (jointly) as a function of $(x,y)\in \Omega\times Y$. Then, the value function $v$ given by \eqref{eq:valuefunction} is Lipschitz and Clarke regular near $\bar{x}$, and $$\partial v(\bar{x}) = \co \{\nabla_x g(\bar{x},\bar{y}): \bar{y}\in S(\bar{x}) \},$$ where $S:\Re^n \rightrightarrows Y$ is given by 
	\begin{equation*}
	S(x) \coloneqq \argmin_{y\in Y} g(x,y).
	\end{equation*}
	\label{thm:danskin}
\end{theorem}

An extension of the above result to the nonsmooth case was also provided in \cite{Bertsekas71} under additional structural assumption on $g$.

\begin{theorem}[Extended Danskin's Theorem]
	Let $g:\Re^n \times Y\to \Re$ be given, where $Y$ is a compact subset of $\Re^m$. Suppose that for a neighborhood $\Omega$ of $\bar{x}\in \Re^n$, $g$ is continuous on $\Omega\times Y$, and for every $y\in Y$, the function $g(\cdot,y)$ is a concave function on $\Re^n$. Then $v$ given by \eqref{eq:valuefunction} is concave on $\Re^n$ and 
	\begin{equation}
	\partial v(\bar{x})= \co \{\xi \in \Re^n: \xi \in \partial_x g(\bar{x},\bar{y})~\text{and}~\bar{y}\in S(\bar{x}) \}.
	\label{eq:P(barx)}
	\end{equation}
	\label{thm:danskin_extended}
\end{theorem}
\ifdefined\submit
\noindent {\textbf{Proof.}}
\else
\begin{proof}
\fi 
	The original result is stated for the maximum value function. Denoting $\ds V(x) \coloneqq \max_{y\in Y} \left( -g(x,y)\right)$, which is a convex function, we have from \cite[Proposition A.22]{Bertsekas71} that $$\partial^c V(\bar{x}) =  \co \{\xi \in \Re^n: \xi \in \partial_x^c (-g(\bar{x},\bar{y}))~\text{and}~\bar{y}\in S(\bar{x}) \}.$$
	Noting that $\partial v = \partial (-V ) = - \partial V= -\partial^c V$, we get the desired result. 
\ifdefined\submit
\hfill \Halmos
\else 
\end{proof}
\fi 

In \cref{thm:danskin_extended_noncompact}, we further extend the above theorem to non-compact sets $Y$ under the assumption that $g$ is a uniform level-bounded function, in the sense defined below; see \cite{RW98}.

\begin{definition}
	\label{defn:uniformlevelbounded}
	A function $h:\Re^m \times \Re^d \to \Re\cup \{+\infty\}$ with values $h(y,z)$ is \emph{level-bounded in $y$ locally uniform in $z$} if for each $z'\in \Re^d$ and $M\in \Re$, there exists an open ball $B$ around $z'$ such that $\{ (y,z) : z\in B, h(y,z)\leq M\}$ is bounded. We also say that $h$ is a \emph{uniformly level-bounded} function.
\end{definition}

We recall the following important result on uniform level-bounded functions; see \cite[Theorem 1.17]{RW98}.

\begin{theorem}%[\citen{Theorem 1.17}{RW98}]
	\label{thm:uniformlevelbounded}
	Let $h:\Re^m\times \Re^d \to \Re \cup \{+\infty \}$ be a proper lower-semicontinuous function that is level-bounded in $y\in \Re^d$ locally uniform in $z\in \Re^d$. Define for all $z\in \Re^d$ the functions
	\[v(z) \coloneqq \min_{y\in \Re^m} h(y,z)\quad \text{and} \quad S(z) \coloneqq \argmin_{y\in \Re^m} h(y,z) \]
	and let $\bar{z}\in \Re^d$. 
	\begin{enumerate}[(a)]
		\item If there exists $\bar{y}\in S(\bar{z})$ such that $h(\bar{y},\cdot)$ is continuous on a set $U$ containing $\bar{z}$, then $v$ is continuous on $U$; and
		\item If $v$ is continuous on a set $U$ containing $\bar{z}$, $\{z^k\}\subseteq U$ such that $z^k\to \bar{z}$ and $\{y^k\}$ is a sequence such that $y^k\in S(z^k)$ for all $k$, then $\{y^k\}$ is bounded and its accumulation points lie \alert{in} $S(\bar{z})$.
	\end{enumerate}
\end{theorem}

\subsection{Smooth approximations}\label{sec:smoothing}
We formally recall some concepts related to smooth approximations. The definitions herein are adopted from \cite{ChenWomersleyYe2011}.
\begin{definition}
	\label{defn:smoothing}
	Let $O\subseteq \Re^d$ be an open set. We say that $\{\phimu:\mu >0\}$ is a \emph{family of smooth approximations} for a continuous function $\phi:O\to \Re$ over the set $O$ if $\phimu:O\to \Re$ is continuously differentiable and if for all $\bar{z}\in O$, we have $ \lim_{z\to \bar{z}, \mu \to 0} \phimu(z) = \phi(\bar{z})$. 
\end{definition}

With the availability of a smooth approximation, a traditional smoothing method is as follows: Given a nonsmooth optimization problem (P) involving a nonsmooth function $g$, we consider a smooth approximate problem (P$_{\mu}$) obtained by replacing $g$ with $\gmu$. This is a convenient approach as there are many powerful algorithms for dealing with smooth optimization problems, which may then be used for solving (P$_{\mu}$). 

Intuitively, we expect that a sequence of solutions (or stationary points) of (P$_{\mu}$) with $\mu$ decreasing to zero will converge to a solution (or stationary point) of the original problem (P). Theoretically, this is achieved through the notion of gradient consistency defined as follows. Note that this property is very important in characterizing stationary points of any smoothing method.

\begin{definition}
	\label{defn:gradconsistency}
	Let $O\subseteq \Re^d$ be an open set and $\{\gmu:\mu >0\}$ be a family of {smoothing functions for a continuous function $g:O\to \Re$ over the set $O$}. This family is said to satisfy the \emph{gradient consistent property} at $\bar{z}\in O$ if 
	\begin{equation*}
	\emptyset \neq \limsup_{z\to\bar{z},\mu\searrow 0} \nabla \gmu (z)  \subseteq \partial g(\bar{z}) , 
	\label{eq:gradientconsistent}
	\end{equation*}
	where 
	\[\limsup_{z\to\bar{z},\mu\searrow 0} \nabla \gmu (z) \coloneqq \left\lbrace \xi\in \Re^d : \exists \{z^k\}, \exists \{\mu_k \} ~\emph{such that}~ z^k\to \bar{z}, \mu_k\searrow 0 ~\emph{and}~ \nabla g_{\mu_k}(z^k) \to \xi \right\rbrace.\] 	
\end{definition}

We now provide some examples of smooth approximations that satisfy gradient consistency. We also highlight the uniform convergence of these smoothing families, to illustrate the satisfiability of the hypotheses of some results in \cref{sec:quadratic,sec:entropy}.

\begin{example} \citet{ChenMangasarian95} proposed a smooth approximation of the plus function $\phi (z) = z_+ = \max\{ z,0\}$, $z\in \Re$, using symmetric density functions $\rho$ with $\kappa \coloneqq \int_{\Re}|s|\rho(s) ds < +\infty$. In particular, they defined
	\[\phimu (z) \coloneqq \int_{\Re} (z-\mu s)_+ \rho (s) ds,\] 
	which satisfies the inequality $0\leq \phimu(z) - \phi(z) \leq \kappa\mu$ for all $z\in\Re$. With this, it is straightforward to see that $\{\phimu : \mu >0 \}$ is a family of smooth approximations of $\phi$ in the sense of \cref{defn:smoothing}, and in addition, $\phi_{\mu} \to \phi$ uniformly on $\Re$. Moreover, it also holds that \alert{$$\limsup_{z\to\bar{z},\mu\searrow 0} \phimu'(z) \begin{cases}
	= 1 & \text{if }\bar{z}> 0\\ 
	= 0 & \text{if }\bar{z}<0\\
	\subseteq [0,1] & \text{if } \bar{z}=0
	\end{cases},$$} and thus gradient consistency is satisfied. This smoothing is particularly useful since many other loss functions and structure-promoting regularizers are expressible in terms of the plus function (see \cite{ANOTC21,ChenXiaojun12}).
	\label{ex:chenmangasarian}
\end{example}

\begin{example}\label{ex:moreau}
	The Moreau envelope of a closed proper function $\phi:\Re^n\to \Re\cup \{ +\infty\}$ is given by 
	\[e_{\mu}\phi (z) \coloneqq \min _{u\in \Re^n} \phi (u) + \frac{1}{2\mu}\norm{z-u}^2.\]
	When $\phi$ is a convex function, $e_{\mu}\phi$ is $(\mu^{-1})-$smooth with $\nabla e_{\mu}\phi(z)=\frac{1}{\mu}(z-p_{\mu}(z))$ where $p_{\mu}(z) \coloneqq \argmin _{u\in \Re^n} \phi (u) + \frac{1}{2\mu}\norm{z-u}^2$ \cite[Theorem 2.26(b)]{RW98}. Moreover, for any open subset $\Omega \subseteq \interior (\dom (\phi))$, we have $0 \leq \phi(z) - e_{\mu}\phi(z) \leq \frac{L_{\phi,\Omega}^2}{2}\mu$ for all $z\in \Omega$, where $L_{\phi,\Omega}$ is the Lipschitz constant of $\phi$ on $\Omega$ (see \cite[Example 9.14]{RW98} and \cite[Theorem 10.51]{Beck17}). Hence, $\{e_{\mu}\phi: \mu >0 \}$ is a family of smooth approximations of $\phi$, and we also note that $e_{\mu}\phi$ converges to $\phi$ uniformly on $\Omega$. Using \cite[Lemma 1]{LiuPongTakeda19} and noting that $\nabla e_{\mu}\phi(z) \in \partial \phi (p_{\mu}(z))$  and $\partial \phi$ is upper semicontinuous, it is not difficult to check that the gradient consistent property is satisfied.
\end{example}

\begin{example}
	Let $g:\Re^n\times \Re^m \to \Re$ be a convex function (that is, convex in $(x,y)$) and level-bounded in $y$ locally uniform in $x$. Then the value function $v$ given by \eqref{eq:valuefunction} is convex \cite[Theorem 10.13]{RW98} and by \cref{ex:moreau}, $\{ \vmu : \mu >0 \}$ where
	\[ \vmu (x) \coloneqq e_{\mu}(v(x)) = \min _{z\in \Re^n} \min_{y\in Y} g(z,y) + \frac{1}{2\mu} \|z-x\|^2\]
	is a family of smooth approximations of the value function satisfying gradient consistency. 
\end{example}

\begin{example}
	Let $g:\Re^n\times \Re^m \to \Re$ be	such that $g(\cdot,y)$ is a concave function for all $y\in Y\subseteq \Re^m$. Then $-v$ is a convex function and we may consider
	\[ \vmu (x) \coloneqq -e_{\mu}(-v(x)) = \max _{z\in \Re^n} \min_{y\in Y} g(z,y) - \frac{1}{2\mu} \|z-x\|^2\]
	as a family of smooth approximations of the value function satisfying gradient consistency. 
\end{example}

Aside from the restrictive structural assumptions on $g$, another main drawback of the smoothing functions in the preceding two examples is that function and gradient evaluations at a single point would require solving a min-min or a max-min optimization problem. 

%\section{Necessary Conditions for Optimality\protect\footnote{Can be put as a subsection of \cref{sec:prelim_nonsmooth}}}
\subsection{Necessary Optimality Conditions}
We recall the definition of stationary point for the value function reformulation \eqref{eq:bilevel_nlp_formulation_epsilon} with $\epsilon\geq 0$. 

\begin{definition}[Stationary point\citen{Definition 4}{YeYuanZengZhang23}]
	Let $(\bar{x},\bar{y})$ be a feasible point of \eqref{eq:bilevel_nlp_formulation_epsilon} with $\epsilon\geq 0$. We say that $(\bar{x},\bar{y})$ is a \emph{stationary point} of \eqref{eq:bilevel_nlp_formulation_epsilon} if there exists $\lambda\geq 0$ such that 
	\begin{equation*}
	\begin{cases}
	0\in \partial f(\bar{x},\bar{y}) + \lambda \partial g(\bar{x},\bar{y}) -  \lambda \partial v(\bar{x})\times \{ 0\} + N_{X\times Y}(\bar{x},\bar{y}) \\
	\lambda (g(\bar{x},\bar{y})-v(\bar{x})-\epsilon )= 0 .
	\end{cases}
	\end{equation*}
	\label{defn:stationary}
\end{definition}
\alert{We also recall the following constraint qualification for \eqref{eq:bilevel_nlp_formulation_epsilon}.}
\alert{\begin{definition}[\citen{Definition 3}{YeYuanZengZhang23}]
	Let $(\bar{x},\bar{y})$ be a feasible point of \eqref{eq:bilevel_nlp_formulation_epsilon} with $\epsilon\geq 0$. We say that the \emph{no nonzero abnormal multiplier constraint qualification (NNAMCQ)} holds at $(\bar{x},\bar{y})$ if either $g(\bar{x},\bar{y})-v(\bar{x})-\epsilon<0$ or $g(\bar{x},\bar{y})-v(\bar{x})-\epsilon=0$ but 
		\begin{equation}
		0 \notin \partial g(\bar{x},\bar{y}) - \partial v(\bar{x})\times \{ 0\} + N_{X\times Y}(\bar{x},\bar{y}).
		\label{eq:nnmamcq}
		\end{equation}
	Let $(\bar{x},\bar{y})\in X\times Y$.  We say that the \emph{extended no nonzero abnormal multiplier constraint qualification (ENNAMCQ)} holds at $(\bar{x},\bar{y})$ if either $g(\bar{x},\bar{y})-v(\bar{x})-\epsilon<0$ or $g(\bar{x},\bar{y})-v(\bar{x})-\epsilon\geq 0$ but \eqref{eq:nnmamcq} holds.
\end{definition}}

\medskip 
\alert{ENNAMCQ is a standard constraint qualification for \eqref{eq:bilevel_nlp_formulation_epsilon} with $\epsilon>0$. It reduces to NNAMCQ when the concerned point belongs to the feasible set of \eqref{eq:bilevel_nlp_formulation_epsilon}. Similar to \cite[Theorem 4]{YeYuanZengZhang23}, we have the following proposition. }
\alert{\begin{proposition}%[\citen{Theorem 4}{YeYuanZengZhang23}]
	\label{prop:localisstationary}
	If $(\bar{x},\bar{y})$ is a local optimal solution to \eqref{eq:bilevel_nlp_formulation_epsilon} with $\epsilon>0$ and NNAMCQ holds at $(\bar{x},\bar{y})$, then $(\bar{x},\bar{y})$ is a stationary point of \eqref{eq:bilevel_nlp_formulation_epsilon}.
\end{proposition}}
	
\alert{\begin{remark}\label{rem:ennamcq}
		\begin{enumerate}[(i)]
			\item Conditions under which ENNAMCQ holds are provided in \cite[Proposition 4.1]{LinXuYe14}. In addition to this, we also see from the proof of \cite[Proposition 8]{YeYuanZengZhang23} that ENNAMCQ is satisfied at $(\bar{x},\bar{y})\in X\times Y$ if 
			\begin{equation}
			\partial g(\bar{x},\bar{y})\subseteq \partial_x g(\bar{x},\bar{y})\times \partial_y g(\bar{x},\bar{y})
			\label{eq:partial_x_partial_y}
			\end{equation}
%			and if $g(\bar{x},\cdot)$ is convex.
			and if $S(\bar{x}) = \{ y \in Y: 	0\in \partial_y g(\bar{x},y)+N_Y(y)  \}$ (for instance,  when $g(\bar{x},\cdot)$ is convex).
			The condition \eqref{eq:partial_x_partial_y} holds when  $g$ is Clarke regular at $(\bar{x},\bar{y})$ by \cite[Proposition 2.3.15]{Clarke83}. In particular, \eqref{eq:partial_x_partial_y} holds when $g$ is convex. 
%			On the other hand, \eqref{eq:0impliessolution} holds when $g(\bar{x},\cdot)$ is convex. 
			\item It is well known that NNAMCQ never holds for \eqref{eq:bilevel_nlp_formulation} \cite[Proposition 7]{YeYuanZengZhang23}. Nevertheless, a result analogous to \cref{prop:localisstationary} can still be established under the so-called ``partial calmness condition''. Since this lies outside the main scope of the present manuscript, we refer the interested reader to \citep{YeZhu95, YeYuanZengZhang23} for further details.
		\end{enumerate}
\end{remark}}

%The following result indicates that local minima of \eqref{eq:bilevel_nlp_formulation_epsilon} are stationary points in the sense of \cref{defn:stationary}. Note that since the nonsmooth MFCQ automatically holds when $\epsilon>0$ (see \cite[Proposition 8]{YeYuanZengZhang23}), the theorem below does not require any CQs, which are typically essential in proving necessary optimality conditions. As mentioned above, lack of a CQ for $\epsilon=0$ is the main hurdle when directly dealing with the reformulation \eqref{eq:bilevel_nlp_formulation}. 

%
%
%We now highlight the importance of gradient consistency when employing smoothing techniques. The following result assumes that stationary points of \eqref{eq:bilevel_nlp_formulation_epsilon_mu} can be obtained. In \cref{sec:quadratic,sec:entropy}, we show that ENNAMCQ holds at points $(\bar{x},\bar{y})\in X\times Y$ for problem \eqref{eq:bilevel_nlp_formulation_epsilon_mu} for particular choices of smooth approximations for the value function. 

\alert{We now highlight the critical role of gradient consistency in the application of smoothing techniques.}

\begin{proposition}
	\label{prop:stationary_smoothing}
	Let $\epsilon>0$. For each $\mu>0$, consider the smooth problem
	\begin{align}
	\begin{array}{rl}
	\ds \min _{(x,y)\in X\times Y}~ & f_{\mu}(x,y) \\
	\text{s.t.} ~& \ds \gmu(x,y) - \vmu (x) \leq \epsilon.
	\end{array}
	\tag{$\text{VFP}_{\epsilon}^{\mu}$}
	\label{eq:bilevel_nlp_formulation_epsilon_mu}
	\end{align}
	where $\{f_{\mu}: \mu>0\}$, $\{\gmu: \mu>0\}$ and $\{\vmu: \mu>0\}$ are families of smooth approximations of $f$, $g$ and $v$, respectively, that all satisfy the gradient consistent property. Let $\{\mu_k\}$ be a sequence of positive numbers with $\mu_k\searrow 0$, and let $(x^k,y^k)$ be a stationary point of \eqref{eq:bilevel_nlp_formulation_epsilon_mu} with $\mu=\mu_k$, and denote by $\lambda_k$ the corresponding Lagrange multiplier. 
	\alert{\begin{enumerate}[(i)]
		\item If $\{(x^k,y^k,\lambda_k)\}$ is bounded, then its accumulation points are stationary points of \eqref{eq:bilevel_nlp_formulation_epsilon}.
		\item 	If $\{(x^k,y^k)\}$ is bounded and $(\bar{x},\bar{y})$ is an accumulation point at which ENNAMCQ for \eqref{eq:bilevel_nlp_formulation_epsilon} holds, then $(\bar{x},\bar{y})$ is a stationary point of \eqref{eq:bilevel_nlp_formulation_epsilon}.
	\end{enumerate}}

\end{proposition}

\ifdefined\submit
\noindent {\textbf{Proof.}}
\else
\begin{proof}
\fi 
	From the KKT conditions of \eqref{eq:bilevel_nlp_formulation_epsilon_mu} with $\mu=\mu_k$, we have that $(x^k,y^k) $ belongs to $ X\times Y$ and satisfies the conditions
	\begin{equation}
	\begin{cases}
	0\in \nabla f_{\mu_k}(x^k,y^k) + \lambda_k \nabla g_{\mu_k} (x^k,y^k) -  \lambda_k \nabla v_{\mu_k}(x^k)\times \{ 0\} +N_{X\times Y} (x^k,y^k) ,\\
	\lambda_k\geq 0, \quad g_{\mu_k}(x^k,y^k)-v_{\mu_k}(x^k)-\epsilon\leq 0, \quad \lambda_k (g_{\mu_k}(x^k,y^k)-v_{\mu_k}(x^k)-\epsilon )= 0,
	\end{cases}
	\label{eq:stationary_mu}
	\end{equation}
	for all $k$. We immediately get the desired result by noting the closedness of $X\times Y$, invoking gradient consistency, and letting $k\to\infty$ in \eqref{eq:stationary_mu}, taking subsequences if necessary. \alert{This proves part (i). To prove part (ii), let $K\subseteq \mathbb{N}$ such that $\{(x^k,y^k)\}_{k\in K}$  converges to $(\bar{x},\bar{y})$. We first claim that $\{\lambda_k\}_{k\in K}$ is bounded. Assume to the contrary that it is unbounded and without loss of generality, $\lambda_k\to +\infty$ as $k\to\infty$ with $k\in K$. Then from \eqref{eq:stationary_mu}, we see that  
		\[0\in \frac{1}{\lambda_k}\nabla f_{\mu_k}(x^k,y^k) +  \nabla g_{\mu_k} (x^k,y^k) -  \nabla v_{\mu_k}(x^k)\times \{ 0\} +N_{X\times Y} (x^k,y^k).\]
	Noting the gradient consistency of hypothesis, taking a subsequence if necessary, and letting $k\to\infty$ with $k\in K$ in the above inclusion, we obtain $0\in \partial g(\bar{x},\bar{y}) - \partial v(\bar{x})\times \{ 0\} + N_{X\times Y}(\bar{x},\bar{y})$. Meanwhile, $\lambda_k>0$ for sufficiently large $k\in K$, and therefore $g_{\mu_k}(x^k,y^k) - v_{\mu_k}(x^k) - \epsilon=0$. Consequently, $g(\bar{x},\bar{y})-v(\bar{x})-\epsilon=0$, yet  $0\in \partial g(\bar{x},\bar{y}) - \partial v(\bar{x})\times \{ 0\} + N_{X\times Y}(\bar{x},\bar{y})$. This is a contradiction since ENNAMCQ holds at $(\bar{x},\bar{y})$ for problem \eqref{eq:bilevel_nlp_formulation_epsilon}. Hence, $\{ \lambda_k\}_{k\in K}$ must be bounded. The result then follows from part (i). This completes the proof. }
\ifdefined\submit 
\hfill \Halmos
\else 
\end{proof}
\fi \alert{\begin{remark}
		In the above proposition, a sequence of KKT points of the smoothed problem \eqref{eq:bilevel_nlp_formulation_epsilon_mu} is shown to have accumulation points that are stationary to \eqref{eq:bilevel_nlp_formulation_epsilon}. By following the same proof as above, it is in fact sufficient to obtain approximate KKT points, in the sense that $(x^k,y^k)\in X\times Y$  and 
		\begin{equation}
		\begin{cases}
		\dist(-\nabla f_{\mu_k}(x^k,y^k) - \lambda_k \nabla g_{\mu_k} (x^k,y^k) +  \lambda_k \nabla v_{\mu_k}(x^k)\times \{ 0\},  N_{X\times Y} (x^k,y^k))\leq \varepsilon_{k,1} ,\\
		\lambda_k\geq -\varepsilon_{k,2}, \quad g_{\mu_k}(x^k,y^k)-v_{\mu_k}(x^k)-\epsilon\leq 0, \quad \lambda_k (g_{\mu_k}(x^k,y^k)-v_{\mu_k}(x^k)-\epsilon )= \varepsilon_{k,3},
		\end{cases}
		\label{eq:stationary_mu_approximate}
		\end{equation}
		where $\varepsilon_{k,i}\to 0$ for $i=1,2,3$. 
\end{remark}}

\medskip 
\alert{Since the above result relies on the ability to compute (approximate) stationary points of \eqref{eq:bilevel_nlp_formulation_epsilon_mu}, it is crucial to ensure that appropriate constraint qualifications hold $-$ particularly because gradient-based algorithms require such conditions to guarantee convergence to stationary points; for instance, see \citep{LinXuYe14,YeYuanZengZhang23}. To this end, in \cref{sec:quadratic,sec:entropy}, we demonstrate that ENNAMCQ holds at points \( (\bar{x}, \bar{y}) \in X \times Y \) for problem \eqref{eq:bilevel_nlp_formulation_epsilon_mu}, under specific choices of smooth approximations to the value function.} \alert{\begin{remark}
		\label{rem:smoothingPG}
As an alternative to computing approximate KKT points of \eqref{eq:bilevel_nlp_formulation_epsilon_mu}, one may consider applying the smoothing projected gradient algorithm proposed in \cite[Algorithm 3.1]{LinXuYe14}, which is designed for nonlinear programs with nonsmooth constraint functions. In that framework, the objective function is assumed to be smooth, while the nonsmooth constraint is replaced by a smooth approximation that satisfies gradient consistency. However, one can show that the same convergence guarantees remain valid even if the objective function is replaced by a gradient-consistent smooth approximation. This reinforces the need for gradient-consistent smooth approximations of $f$, $g$ and $v$.
\end{remark} }

\medskip 
\alert{Motivated by \cref{prop:stationary_smoothing,rem:smoothingPG}, we propose smoothing techniques for the value function that satisfy gradient consistency in \cref{sec:quadratic,sec:entropy}.}

\section{Quadratic Regularization}
\label{sec:quadratic}
In this section, we consider the case when $g(x,\cdot)$ is a convex function for each $x\in O_X$, which will be our standing assumption. We also assume that we have a family $\{\gmu:\mu>0\}$ of smooth approximations of $g$. In this case, note that $\{\gtildemu: \mu > 0\}$ where
\begin{equation}
\gtildemu (x,y) \coloneqq \gmu (x,y) + \frac{\mu}{2}\|y\|^2,
\label{eq:smoothing_stronglyconvex}
\end{equation}
is likewise a family of smoothing functions for $g$. 

\subsection{Smooth approximation of the value function}
Intuitively, we see that we may approximate the value function $v$ by replacing the corresponding objective function $g$ by its smooth approximation $\gtildemu$ (or $\gmu$). Here, we consider the functions $\alert{\vmuquad}:\Re^n\to \Re$ and $\smuquad: \Re^n \rightrightarrows \Re^m$ given by
\begin{equation}
\alert{\vmuquad}(x) \coloneqq \min_{y \in Y} \gtildemu (x,y) \quad \text{and} \quad \smuquad(x) \coloneqq \argmin_{y\in Y} \gtildemu (x,y) 
\label{eq:valuefunction_mu}
\end{equation}
where $\gtildemu$ is given by \eqref{eq:smoothing_stronglyconvex}. In the following discussion, we show that $\{ \alert{\vmuquad} : \mu>0\}$ is indeed a family of smoothing functions for $v$. \alert{Assumptions to guarantee differentiability of $\vmuquad$ are summarized below.}
\alert{\begin{assumption}\label{assume:quadreg}
	The following holds:
	\begin{enumerate}[(a)]
		\item For any open set $\Omega$ with $\cl (\Omega)\subseteq O_X$, $\gmu(x,y)$ is uniformly bounded below over $\cl (\Omega)\times Y$ for any $\mu>0$. 
		\item $\gmu(x,\cdot)$ is convex for any $x\in O_X$ and $\mu>0$.
	\end{enumerate}
\end{assumption}}

\begin{remark}\label{remark:uniformbounded} 
	We make some comments about the uniform boundedness assumption in \alert{\cref{assume:quadreg}(a)}.
	\begin{enumerate}
		\item If $Y$ is compact, this assumption automatically holds by the continuity of $\gmu$. 
		\item 	Note that the uniform boundedness assumption holds for a wide class of bilevel problems with unbounded $Y$. For instance, consider the class of hyperparameter problems in \cite{ANOTC21}, where the lower-level function takes the form 
		\begin{equation}
		g(x,y) = \ell (y)+ \sum_{i=1}^{n}x_i p_i(y),
		\label{ex:hyperparameter_learning}
		\end{equation}
		where $\ell$ corresponds to a nonnegative loss function, $x $ is a hyperparameter, and $p_i$ are (possibly nonsmooth) regularizers which are nonnegative. The function $g$ along with the family of smoothing functions derived in \cite{ANOTC21} are uniformly bounded below by $0$ over $\Re^n_+ \times \Re^m$. Hence, uniform boundedness holds when $Y=\Re^m$ and $X=[\varepsilon_1,\infty)\times \cdots \times [\varepsilon_n,\infty)$, where $\varepsilon_i>0$ for all $i$. 
	\end{enumerate}
\end{remark}

\begin{remark}
	\label{remark:moreau}
	If $g(x,\cdot)$ is convex for each $x\in O_X$, then its Moreau envelope $e_{\mu}g(x,y)$ is convex in $y$ for any $x\in O_X$, which satisfies the convexity requirement of the smoothing function $\gmu$ in \alert{\cref{assume:quadreg}(b)}. 
\end{remark}

We mention that the differentiability of the value function associated with a quadratically regularized \textit{smooth function} of $(x,y)$ that is convex in $y$ has already been noted in other works; for instance, see \citep{Liu23}. Since this is precisely the structure of the function $\gtildemu$ (i.e., it is the sum of $\gmu$, a smooth function convex in $y$, and the quadratic term $\mu\|y\|^2/2$), then $\alert{\vmuquad}$ given in \eqref{eq:valuefunction_mu} must be smooth. \alert{We summarize this fact in the following proposition for completeness.} Moreover, we emphasize as well that in this paper, we prove that the obtained function $\alert{\vmuquad}$ is a smooth approximation of $v$ in the sense of \cref{defn:smoothing}, which, to the best of our knowledge, has not been shown in the literature. Even more importantly, we establish the gradient consistency of this smoothing family, which is instrumental in developing smoothing algorithms that converge to stationary points of the bilevel problem (see \cref{prop:stationary_smoothing}).

\begin{proposition}
	\label{prop:vmu_C1}
	\alert{Suppose that \cref{assume:quadreg} holds}. Then $\alert{\vmuquad}$ is continuously differentiable on $O_X$ with $\nabla \alert{\vmuquad} (x) = \nabla_x \gmu (x,\smuquad (x))$. 
\end{proposition}
\ifdefined\submit
\noindent {\textbf{Proof.}}
\else
\begin{proof}
\fi 
\alert{This follows directly from \cite[Proposition 1]{Liu23}.}
\ifdefined\submit
\hfill \Halmos
\else 
\end{proof}
\fi

\alert{Having established the differentiability of $ \vmuquad $, we next verify that any point $ (\bar{x}, \bar{y}) \in X \times Y $ satisfies ENNAMCQ for \eqref{eq:bilevel_nlp_formulation_epsilon_mu} with $ v_{\mu} = \vmuquad $.}
\alert{\begin{proposition}\label{prop:ennamcq_quad}
Suppose that \cref{assume:quadreg} holds, and let $v_{\mu}=\vmuquad$ in \eqref{eq:bilevel_nlp_formulation_epsilon_mu} with $\epsilon>0$. Then ENNAMCQ holds at any point $(\bar{x},\bar{y})\in X\times Y$, that is, either $g_{\mu}(\bar{x},\bar{y})-\vmuquad (\bar{x})<\epsilon$ or $g_{\mu}(\bar{x},\bar{y})-\vmuquad (\bar{x})\geq \epsilon$ but $0\notin \nabla g_{\mu}(\bar{x},\bar{y}) - \nabla \vmuquad (\bar{x})\times \{0\} + N_{X\times Y}(\bar{x},\bar{y}). $
\end{proposition}}
\alert{\ifdefined\submit
\noindent {\textbf{Proof.}}
\else
\begin{proof}
	\fi
Suppose $g_{\mu}(\bar{x},\bar{y})-\vmuquad (\bar{x})\geq \epsilon$, and suppose, for the sake of contradiction, that $0\in \nabla g_{\mu}(\bar{x},\bar{y}) - \nabla \vmuquad (\bar{x})\times \{0\} + N_{X\times Y}(\bar{x},\bar{y})$. From the latter condition, we have $0\in \nabla_y g_{\mu}(\bar{x},\bar{y}) + N_Y(\bar{y})$, which by \cref{assume:quadreg}(b) implies that $\bar{y}\in \argmin_{y\in Y} g_{\mu}(\bar{x},y)$. It follows that 
\[\vmuquad (\bar{x}) = g_{\mu}(\bar{x},S_{\mu}(\bar{x})) + \frac{\mu}{2}\norm{S_{\mu}(\bar{x})}^2 \geq g_{\mu}(\bar{x},\bar{y}) + \frac{\mu}{2}\norm{S_{\mu}(\bar{x})}^2. \]
Hence, 
\[0\geq g_{\mu}(\bar{x},\bar{y}) - \vmuquad (\bar{x}) + \frac{\mu}{2}\norm{S_{\mu}(\bar{x})}^2 \geq \epsilon + \frac{\mu}{2}\norm{S_{\mu}(\bar{x})}^2, \]
and so $\epsilon=0$, which is a contradiction. Therefore, $0\notin \nabla g_{\mu}(\bar{x},\bar{y}) - \nabla \vmuquad (\bar{x})\times \{0\} + N_{X\times Y}(\bar{x},\bar{y})$.
	\ifdefined\submit
	\hfill \Halmos
	\else 
\end{proof}
\fi }

\alert{We now show that $\vmuquad$ is a smooth approximation of $v$ in the sense of \cref{defn:smoothing}. To this end, we need some additional assumptions.}
\alert{\begin{assumption}\label{assume:quadred_2}
	$\gmu(x,y)$ is continuous as a function of $(x,y,\mu)\in \Re^n\times \Re^m\times \Re_{++}$.
\end{assumption}
\begin{assumption}\label{assume:quadred_3}
	$g$ is level-bounded in $y$ and locally uniform in $x$.
\end{assumption}
\begin{assumption}\label{assume:quadred_4} For any open set $\Omega$ with $\cl (\Omega)\subseteq O_X$, $\gmu\to g$ uniformly on $\cl (\Omega)\times Y$.
\end{assumption}}
\alert{Examples of smoothing approximations that satisfy \cref{assume:quadred_2,assume:quadred_4} are given in \cref{ex:chenmangasarian,ex:moreau}. On the other hand, \cref{assume:quadred_3} is a standard condition when dealing with non-compact constraint sets (for instance, see \cite{Liu23}).}

\begin{theorem}
	\label{thm:smoothingfunction_valuefunction_noncompact}
	\alert{Suppose that \cref{assume:quadreg,assume:quadred_2,assume:quadred_3,assume:quadred_4} hold.} Then $\{ \alert{\vmuquad} : \mu > 0\}$ is a family of smoothing functions over $O_X$ for the value function $v$. 
\end{theorem}
\ifdefined\submit
\noindent {\textbf{Proof.}}
\else
\begin{proof}
\fi
	Since $\alert{\vmuquad}$ is continuously differentiable on $O_X$ by \cref{prop:vmu_C1}, it is enough to show that $\lim_{x\to \bar{x}, \mu\to 0} \alert{\vmuquad}(x) = v(\bar{x})$ for any $\bar{x}\in O_X$. To this end, define 	$h:\Re^m\times \Re^n\times \Re\to \Re\cup\{+\infty \}$ as
	\begin{equation}
	h(y;x,\mu) \coloneqq \begin{cases}
	g_{|\mu|} (x,y) + \frac{|\mu|}{2}\|y\|^2 + \delta_{\cl (\Omega) \times Y}(x,y) & \text{if}~\mu\neq 0 \\
	g(x,y) + \delta_{\cl (\Omega)\times Y} (x,y) & \text{if}~\mu=0
	\end{cases},
	\label{eq:h_extended}
	\end{equation} 
	where $\Omega$ is an open set containing $\bar{x}$ with $\cl (\Omega)\subseteq O_X$. 
	%	For all $(x,\mu) \in \cl (B_{\delta}(\bar{x}))\times \Re$, define
	For all $(x,\mu) \in \cl (O_X)\times \Re$, define
	\begin{equation}
	\alert{\vmuquad} (x) \coloneqq \min_{y\in \Re^m} 	h(y;x,\mu) \quad \text{and} \quad \smuquad (x) \coloneqq \argmin_{y\in \Re^m} 	h(y;x,\mu),
	\label{eq:v_intermsof_h}
	\end{equation}
	which agree on the set $\cl (\Omega)$ with the definitions of $\alert{\vmuquad}$ and $\smuquad$ given in \eqref{eq:valuefunction_mu} when $\mu>0$. We now claim that $\alert{\vmuquad} (x)$ viewed as a function of $(x,\mu)\in \Re^n\times \Re$ is continuous on a neighborhood of $(\bar{x},0)$. To this end, \cref{thm:uniformlevelbounded} asserts that it is enough to show that (i) $h$ given by \eqref{eq:h_extended} is a proper lower semicontinuous function; (ii) $h$ is level-bounded in $y$ locally uniform in $(x,\mu)$; and (iii) there exists $\bar{y}\in S_{0}(\bar{x}) = S(\bar{x})$ such that $h(\bar{y};x,\mu)$ is a continuous function of $(x,\mu)$ in some neighborhood of $(\bar{x},0)$.

	To prove (i), \alert{we first note from \cref{assume:quadred_2} that $g_{\mu}(x,y)$ is continuous as a function of $(x,y,\mu)$ on $\Re^n\times \Re^m \times \Re_{++}$ and so} $h$ is lower semicontinuous at any $(y,x,\mu)$ with $\mu \neq 0$. On the other hand, since $g_{\mu}$ is a smoothing function of $g$, it follows from \cref{defn:smoothing} that $h$ is also lower semicontinuous at any point $(y,x,0)$ with $(y,x)\in \Re^m\times \Re^n$. In fact, we have that $h$ is continuous on the set $Y\times \cl (O_X)\times \Re$. From this, we also obtain claim (iii) by noting that $S_0(\bar{x}) = S(\bar{x})\neq \emptyset$ due to our uniform level-boundedness assumption on $g$ \alert{(\cref{assume:quadred_3})}. It remains to show that (ii) holds. We need to show that given any $M\in \Re$, there exists $\delta>0$ such that $L(M;\delta)\coloneqq \{(y,x,\mu): x\in B_{\delta}(\bar{x}), |\mu|<\delta , ~h(y;x,\mu)\leq M\}$ is bounded.  By the uniform convergence \alert{of $\gmu$ as given in \cref{assume:quadred_4}}, there exists $\delta_1>0$ such that $\gmu (x,y)>g(x,y)-1$ for all $(x,y)\in \cl(\Omega)\times Y$ and $\mu\in (0,\delta_1)$.  Then $h(y;x,\mu) \geq g(x,y) - 1$ for all $(x,y,\mu)\in  \Re^n\times \Re^m\times (-\delta_1,\delta_1)$. Since $g$ is level bounded in $y$ locally uniform in $x$, there exists $\delta_2>0$ such that $\{ (y,x): x\in B_{\delta_2}(\bar{x}), ~g(x,y)\leq M+1 \}$ is bounded. Choosing $\delta\coloneqq \min\{\delta_1,\delta_2\}$, the resulting set $L(M;\delta)$ is in turn bounded, and therefore $h$ is level bounded in $y$ locally uniform in $(x,\mu)$.   Hence, we have completed the proof that $\alert{\vmuquad}(x)$ is continuous at $(\bar{x},0)$ for any $\bar{x}\in O_X$. Consequently, by definition of continuity, 
	\[ \lim_{x\to \bar{x}, \mu\to 0} v_{\mu}(x) = v_0 (\bar{x}) = \min_{y \in \Re^m} ~h(y;\bar{x},0) = \min_{y\in \Re^m} g(\bar{x},y)+\delta_Y (y) = \min_{y\in Y} g(\bar{x},y) = v(\bar{x}), \]
	showing that $\alert{\vmuquad}$ is indeed a smoothing function for \alert{$v$}. 
\ifdefined\submit
\hfill \Halmos
\else 
\end{proof}
\fi

\subsection{Danskin-type theorems and gradient consistency}
We have shown above that regularized smoothing approximations $\{\gtildemu : \mu>0 \}$ of $g$ lead to smoothing functions $\{\alert{\vmuquad} : \mu>0 \}$ for the value function $v$. Unfortunately, the gradient consistent property is not necessarily inherited by the approximations for $v$; that is,  $\{\alert{\vmuquad} : \mu>0 \}$ does not necessarily possess gradient consistent property when $\{\gtildemu : \mu>0 \}$ has this property. Fortunately, we can obtain gradient consistency for several classes of functions $g$. 

\subsubsection*{Danskin-type theorems}
To prove gradient consistency of  $\{\alert{\vmuquad} : \mu>0 \}$, it is important to determine what exactly the elements of the Clarke subdifferential\footnote{Note that $v$ is Lipschitz continuous so that its Clarke subdifferential is well-defined. This fact follows from the Lipschitz continuity of $g$. Indeed, given $x$ and $\bar{x}$, denote $y\in S(x) $ and $\bar{y}\in S(\bar{x})$. Then $v(x)-v(\bar{x}) = g(x,y)-g(\bar{x},\bar{y}) \leq g(x,\bar{y})-g(\bar{x},\bar{y}) \leq L \|x-\bar{x}\|$, where $L$ is the Lipschitz constant of $g$. Reversing the roles of $x$ and $\bar{x}$, we immediately get $|v(x)-v(\bar{x})|\leq L\|x-\bar{x}\|$.} of $v$ are. For compact sets $Y$, the subdifferential $\partial v$ is characterized by Danskin's Theorem (\cref{thm:danskin}) for smooth functions $g$, and its extension \cref{thm:danskin_extended} for a nonsmooth $g$ that is concave in $x$. We now establish extensions of these results for non-compact sets $Y$ by assuming uniform level-boundedness of $g$. 

We start with the following extension of the classical Danskin's Theorem (\cref{thm:danskin}). This result should be well-known, and we provide a simple proof for self-containednenss.

\begin{theorem}
	Let $g:\Re^n \times Y\to \Re$ be given, where $Y$ is a closed subset of $\Re^m$ and $g$ is level-bounded in $y$ locally uniform in $x$. Suppose that for a neighborhood $\Omega$ of $\bar{x}\in \Re^n$, the derivative $\nabla_x g(x,y)$ exists and is continuous (jointly) as a function of $(x,y)\in \Omega\times Y$. Then the value function $v$ given by \eqref{eq:valuefunction} is Lipschitz and Clarke regular near $\bar{x}$, and $$\partial v(\bar{x}) = \co \{\nabla_x g(\bar{x},\bar{y}): \bar{y}\in S(\bar{x}) \},$$ where $	S(x) \coloneqq \argmin_{y\in Y} g(x,y).$
	\label{thm:danskin_noncompact}
\end{theorem}
\ifdefined\submit
\noindent {\textbf{Proof.}}
\else
\begin{proof}
\fi
	Set $M\coloneqq  \max_{x\in \cl (B_{1}(\bar{x}))} g(x,y')$ for some fixed and arbitrarily chosen $y'\in Y$. Note that by continuity of $g$, $M$ is finite. By uniform level-boundedness of $g$, there exists some ball $B_{\alpha}(\bar{x})$ with $\alpha <1$ such that $\{ (x,y) : x\in \cl (B_{\alpha}(\bar{x})), g(x,y) \leq M\}$ is bounded. Since $\alpha<1$, then
	\[\bigcup_{x\in  \cl (B_{\alpha}(\bar{x}))} S(x) \subseteq \bigcup_{x\in \cl (B_{\alpha}(\bar{x}))} \{ y\in Y : g(x,y) \leq M\}. \]
	It follows that $\hat{Y} \coloneqq  \bigcup_{x\in  \cl (B_{\alpha}(\bar{x}))} S(x)$ is a bounded set. If we have a sequence $\{y^k\}\subseteq \hat{Y}$ such that $y^k \to y^*$, then there exists $x^k\in \cl (B_{\alpha}(\bar{x}))$ such that $g(x^k,y^k)\leq g(x^k,y)$ for all $k$ and for all $y\in Y$. We may assume without loss of generality that $x^k\to x^*$ where $x^*\in  \cl (B_{\alpha}(\bar{x}))$. Noting the continuity of $g$, we have $g(x^*,y^*)\leq g(x^*,y)$ for all $y\in Y$ so that $y^*\in S(x^*)$. Hence, $y^*\in \hat{Y}$ and therefore $\hat{Y}$ is in fact compact. Consider now the functions 
	\[ \hat{v}(x) \coloneqq \min_{y\in \hat{Y}} g(x,y) \quad \text{and} \quad \hat{S}(x) \coloneqq \argmin_{y\in \hat{Y}} g(x,y) .\]
	Since the global minimizers of $g(x,\cdot)$ are contained in $\hat{Y}$ for any $x\in \cl (B_{\alpha}(\bar{x}))$, it follows that $\hat{S}\equiv S$ on $\cl (B_{\alpha}(\bar{x}))$. In particular, $\hat{v}\equiv v$ on the open set $B_{\alpha}(\bar{x})$ so that their Clarke subdifferentials agree. Hence, together with Danskin's theorem, $\partial v(\bar{x})  = \partial \hat{v}(\bar{x}) = \co \{ \nabla_x g(\bar{x},\bar{y}) : \bar{y}\in \hat{S}(\bar{x}) \}= \co \{ \nabla_x g(\bar{x},\bar{y}) : \bar{y}\in {S}(\bar{x}) \}$. This completes the proof. 
\ifdefined\submit
\hfill \Halmos
\else 
\end{proof}
\fi 

On the other hand, the following theorem is an extension of \cref{thm:danskin_extended} that is applicable for unbounded closed sets $Y$ provided that $g$ is uniformly level-bounded. Moreover, we also relax the concavity assumption to weak concavity. Recall that given a convex set $\Omega$, a function $h$ is $\rho$-weakly convex on $\Omega$ if $h+\frac{\rho}{2}\|\cdot \|^2$ is convex on $\Omega$, where $\rho >0$. We say that $h$ is $\rho$-weakly concave on $\Omega$ if $-h$ is a weakly convex function on $\Omega$. 

\begin{theorem}
	\label{thm:danskin_extended_noncompact}
	Let $g:\Re^n \times Y\to \Re$ be level-bounded in $y$ locally uniform in $x$, where $Y$ is a closed subset of $\Re^m$. Suppose that for a convex neighborhood $\Omega$ of $\bar{x}\in \Re^n$ and an open set $O_Y$ containing $Y$, $g$ is Lipschitz continuous on $\Omega\times O_Y$ and $g(\cdot,y)$ is a $\rho$-weakly concave function on $\Omega$ for every $y\in O_Y$. Then $v$ given by \eqref{eq:valuefunction} is weakly concave on $\Omega$ and 
	\begin{eqnarray}
	\partial v(\bar{x}) = \co \{\xi \in \Re^n: \xi \in \partial_x g(\bar{x},\bar{y})~\text{and}~\bar{y}\in S(\bar{x}) \}.
	\label{eq:P(barx)_noncompact}
	\end{eqnarray}
\end{theorem}

\ifdefined\submit
\noindent {\textbf{Proof.}}
\else
\begin{proof}
\fi
	Define $h\coloneqq -g$, so that $V(x) \coloneqq \max_{y\in Y} h(x,y) = - v(x)$ and $h(\cdot,y)$ is $\rho$-weakly convex on $\Omega$ for each $y\in O_Y$. It follows that $\hat{h}(\cdot,y) \coloneqq h(\cdot,y)+ \frac{\rho}{2}\|\cdot \|^2$ is a convex function on $\Omega$, and therefore the maximum over all $y\in Y$ is convex on $\Omega$. That is, 
	\[\max_{y\in Y} \hat{h}(x,y) = \max_{y\in Y} \left(h(x,y)+ \frac{\rho}{2}\|x\|^2\right) = \left( \max_{y\in Y} h(x,y) \right) + \frac{\rho}{2} \|x\|^2 = V(x) + \frac{\rho}{2} \|x\|^2 \eqqcolon \hat{V}(x)\] 
	is convex on $\Omega$, which proves that $V$ is $\rho$-weakly convex on $\Omega$. Hence, $v$ is $\rho$-weakly concave on $\Omega$.  
	
	Let $d\in \Re^n$, $t_k\searrow 0$, and denote $x^k \coloneqq \bar{x} + t_k d$ and $y^k\in S(x^k)$. Note that given any $y\in S(\bar{x})$, 
	%	\begin{align*}
	%	\frac{\hat{V}(\bar{x} + t_k d)-\hat{V}(\bar{x})}{t_k} & = \frac{h(\bar{x}+t_k d,y^k)+\frac{\rho}{2}\norm{\bar{x}+t_kd}^2-h(\bar{x},y)-\frac{\rho}{2}\|\bar{x}\|^2}{t_k} \\ 
	%	& \geq \frac{h(\bar{x}+t_k d,y)+\frac{\rho}{2}\norm{\bar{x}+t_kd}^2-h(\bar{x},y)-\frac{\rho}{2}\|\bar{x}\|^2}{t_k} \\
	%	& =  \frac{h(\bar{x}+t_k d,y)-h(\bar{x},y)}{t_k} + \frac{\rho}{2t_k}(2t_k\bar{x}^\top d +t_k^2\|d\|^2) 
	%	\end{align*}
	\begin{align*}
	\frac{\hat{V}(\bar{x} + t_k d)-\hat{V}(\bar{x})}{t_k} & = \frac{\hat{h}(\bar{x}+t_k d,y^k)-\hat{h}(\bar{x},y)}{t_k}  \geq \frac{\hat{h}(\bar{x}+t_k d,y)-\hat{h}(\bar{x},y)}{t_k} 
	\end{align*}
	Letting $k\to \infty$, we have
	\begin{equation}
	\hat{V}'(\bar{x};d) \geq \hat{h}_{y} '(\bar{x};d) , \quad \forall y\in S(\bar{x}),
	\label{eq:Vhat>=}
	\end{equation}
	where $\hat{h}_{y} '(\bar{x};d)$ denotes the directional derivative of $h(\cdot,y)$ at $\bar{x}$ in the direction $d$. Meanwhile, we also have
	\begin{align}
	\frac{\hat{V}(\bar{x} + t_k d)-\hat{V}(\bar{x})}{t_k} & = \frac{\hat{h}(\bar{x}+t_k d,y^k)-\hat{h}(\bar{x},y)}{t_k} \leq\frac{\hat{h}(\bar{x}+t_k d,y^k)-\hat{h}(\bar{x},y^k)}{t_k}.
	\label{eq:Vhat<=}
	\end{align}
	Note that since $g$ is level-bounded in $y$ locally uniform in $x$, we have from \cref{thm:uniformlevelbounded} that $\{y^k\}$ is bounded with accumulation points lying in $S(\bar{x})$. Hence, we may assume without loss of generality that $y^k\to \bar{y}\in S(\bar{x})$. Letting $k\to\infty$ in \eqref{eq:Vhat<=} and from the definition of the Clark generalized directional derivative, we see that 
	\begin{equation}
	\hat{V}'(\bar{x};d)\leq \hat{h}^{\circ} (\bar{x},\bar{y}; d,0).
	\label{eq:Vhat<=2}
	\end{equation}
	From the max formula \eqref{eq:maxformula} for Clarke generalized directional derivative, we have
	\begin{equation}
	\hat{h}^{\circ} (\bar{x},\bar{y}; d,0) = \max \{ \lla \zeta, (d,0)\rla : \zeta \in \partial \hat{h}(\bar{x},\bar{y}) \} = \max \{ \lla \zeta_x , d\rla : \zeta_x \in \pi_x (\partial \hat{h}(\bar{x},\bar{y}))\} ,
	\label{eq:h_hat_circ}
	\end{equation}
	where $\pi_x$ is defined as in \eqref{eq:pi_x}.
	%	where $A(\bar{x},\bar{y}) \coloneqq \{ \zeta_x : \exists \zeta_y~\text{s.t. }(\zeta_x,\zeta_y)\in \partial \hat{h}(\bar{x},\bar{y})\}$.
	On the other hand, by the convexity of $\hat{h}(\cdot,y)$ for all $y$ on some neighborhood of $\bar{y}$, we have from \cite[Proposition 2.5.3]{Clarke83} that $\pi_x (\partial \hat{h}(\bar{x},\bar{y})) \subseteq \partial_x \hat{h}(\bar{x},\bar{y})$. From \eqref{eq:h_hat_circ}, we conclude that 
	\begin{equation}
	\hat{h}^{\circ} (\bar{x},\bar{y}; d,0) \leq \max \{ \lla \zeta_x, d\rla : \zeta_x \in \partial_x \hat{h}(\bar{x},\bar{y})\} = \hat{h}_{\bar{y}}'(\bar{x};d),
	\label{eq:h_hat_circ_leq}
	\end{equation}
	where the last equation holds by the max formula for the (usual) directional directive for the convex function $\hat{h}_{\bar{y}}$.  Combining \eqref{eq:Vhat>=}, \eqref{eq:Vhat<=2} and \eqref{eq:h_hat_circ_leq} while noting that $\bar{y}\in S(\bar{x})$, we conclude that 
	\begin{align*}
	\hat{V}'(\bar{x};d) = \max_{y\in S(\bar{x})} \hat{h}_y' (\bar{x};d) & = \max_{y\in S(\bar{x})} \max_{\zeta\in \partial_x \hat{h}_y (\bar{x})} \lla \zeta , d\rla  \\ 
	& = \max \left\lbrace \lla \zeta, d\rla : \zeta \in \bigcup_{y\in S(\bar{x})}\partial_x \hat{h}(\bar{x},y) \right\rbrace \\
	& = \max \left\lbrace \lla \zeta, d\rla : \zeta \in \cl \left( \co \left( \bigcup_{y\in S(\bar{x})}\partial_x \hat{h}(\bar{x},y) \right) \right) \right\rbrace ,
	\end{align*}
	where the last equation holds by \cite[Lemma 2.35]{Beck17}. Noting the  upper semicontinuity of $\partial h $ and compactness of $S(\bar{x})$, we obtain the compactness of $\bigcup_{y\in S(\bar{x})}\partial_x \hat{h}(\bar{x},y)$, and therefore its convex hull is likewise compact \cite[Theorem 3.20]{Rudin91}. Thus, we may equivalently write 
	\[ \hat{V}'(\bar{x};d)  = \max \left\lbrace \lla \zeta, d\rla : \zeta \in  \co \left( \bigcup_{y\in S(\bar{x})}\partial_x \hat{h}(\bar{x},y) \right)\right\rbrace. \]
	Meanwhile, we have from the convexity of $\hat{V}$ that $ \hat{V}'(\bar{x};d) = \max_{\xi \in \partial \hat{V}(\bar{x})} \lla \xi, d \rla $. By the closedness and convexity of both $\partial \hat{V}(\bar{x})$ and $\co \left( \bigcup_{y\in S(\bar{x})}\partial_x \hat{h}(\bar{x},y) \right)$, we have from \cite[Lemma 2.34]{Beck17} that 
	\[ \partial \hat{V}(\bar{x}) = \co \left( \bigcup_{y\in S(\bar{x})}\partial_x \hat{h}(\bar{x},y) \right).\]
	Subtracting $\rho \bar{x}$ from both sides of the above equation, we see that $\partial V(\bar{x}) =\co \left( \bigcup_{y\in S(\bar{x})}\partial_x h(\bar{x},y) \right)$. From here, we immediately obtain \eqref{eq:P(barx)_noncompact}, thus completing the proof. 
\ifdefined\submit
\hfill \Halmos
\else 
\end{proof}
\fi

We also take note of the following special case. %which should also be well-known, but a proof is provided for completeness. 

\begin{theorem}
	\label{thm:danskin_convex}
	Let $g(x,y)$ be convex in $(x,y)$ and level-bounded in $y$ locally uniform in $x$. If the partial differentiation formula 
	\begin{equation}
	\partial g(\bar{x},\bar{y}) = \partial _x g(\bar{x},\bar{y}) \times \partial _yg(\bar{x},\bar{y}) 
	\label{eq:partialformula}
	\end{equation}
	holds for all $\bar{y}\in S(\bar{x})$, then
	\begin{equation*}
	\partial v(\bar{x}) = \partial_x g(\bar{x},\bar{y}),
	\end{equation*}
	where $\bar{y}$ can be arbitrarily chosen from $S(\bar{x})$. 
\end{theorem}
\ifdefined\submit
\noindent {\textbf{Proof.}}
\else
\begin{proof}
\fi
%	Since $g$ is convex in $(x,y)$, we have from \cite[Theorem 10.13]{RW98} that $\partial v(\bar{x}) = \{ \xi \in \Re^n : (\xi,0)\in \partial \left( g(\bar{x},\bar{y}) + \delta_Y (\bar{y}) \right) \}$ where $\bar{y}$ can be arbitrarily chosen from $S(\bar{x})$. Meanwhile, $\partial \left( g(\bar{x},\bar{y}) + \delta_Y (\bar{y}) \right) = \partial g(\bar{x},\bar{y}) + \left( \{0\}\times N_Y(\bar{y}) \right) $ by \cite[Theorem 3.36(a)]{Beck17} and \cite[Corollary 1 of Theorem 2.9.8]{Clarke83}, noting that $g$ is Lipschitz continuous. Further applying \eqref{eq:partialformula}, we then obtain
%	\begin{equation*}
%	\partial v(\bar{x}) = \left\lbrace \xi\in \Re^n : (\xi,0) \in \partial_xg (\bar{x},\bar{y}) \times \left( \partial_y g (\bar{x},\bar{y}) +N_Y(\bar{y}) \right) \right\rbrace = \partial_x g(\bar{x},\bar{y})\quad \forall \bar{y}\in S(\bar{x}),
%	%\label{eq:valuefunction_convexobjective}	
%	\end{equation*}
%	where the last equality holds by the optimality of $\bar{y}$, noting the convexity of $g$ and $Y$.	
\alert{This directly follows from \cite[Theorem 3]{YeYuanZengZhang23}.}
\ifdefined\submit
\hfill \Halmos
\else 
\end{proof}
\fi 

As noted in \cref{rem:ennamcq}, we always have  $	\partial g(\bar{x},\bar{y}) \subseteq \partial _x g(\bar{x},\bar{y}) \times \partial _yg(\bar{x},\bar{y}) $ for a convex function $g$. 
%More generally, this inclusion holds when $g$ is Clarke regular at $(\bar{x},\bar{y})$ by \cite[Proposition 2.3.15]{Clarke83}. 
Conditions under which equality holds, as required by condition \eqref{eq:partialformula}, are provided in \cite[Proposition 1]{YeYuanZengZhang23}.

\subsubsection*{Gradient consistency}
With the above characterizations of the Clarke subdifferential of the value function, we can now derive the gradient consistency of the proposed smoothing family $\{ \alert{\vmuquad} : \mu>0\}$ given in \eqref{eq:valuefunction_mu}. Our main tool is the following simple lemma. 
%\begin{lemma}
%	\label{lemma:s_mu}
%	Suppose that the hypotheses of \cref{thm:smoothingfunction_valuefunction_noncompact} hold. Then for any $\bar{x}\in O_X$, the sequence $\{ S_{\mu_k} (x^k)\}$ where $\mu_k\searrow 0$ and $x^k\to \bar{x}$ is bounded and its accumulation points belong to $S(\bar{x})$. Consequently, 
%	\[\lim_{x\to \bar{x}, \mu\searrow 0} \dist (\smuquad (x), S(\bar{x})) = 0 \quad \forall \bar{x}\in O_X.\]
%\end{lemma}
%\begin{proof}
%	This directly follows from 
%	\cref{thm:smoothingfunction_valuefunction_noncompact} and \cref{thm:uniformlevelbounded}(b).
%\end{proof}

\begin{lemma}
	\label{thm:gradientconsistentproperty}
	\alert{Suppose that \cref{assume:quadreg,assume:quadred_2,assume:quadred_3,assume:quadred_4} hold.} and that $\{\gmu: \mu>0\}$ is a smoothing family of functions that satisfies the gradient consistent property. Let $\bar{x}\in O_X$ and suppose that 
	\begin{equation}
	\pi_x (\partial g\left(\bar{x},\bar{y})\right) \subseteq \partial_x g(\bar{x},\bar{y}) \quad \forall \bar{y}\in S(\bar{x}),
	\label{eq:projection_first_element_subset}
	\end{equation}
	and $\{\nabla \gmu (x,y) : (x,y,\mu) \in \Theta\}$ is bounded\footnote{This condition is easily satisfied by several smoothing functions; see \cref{sec:smoothing} and \cite{ANOTC21,ChenXiaojun12}.} for some neighborhood $\Theta$ of $(\bar{x},\bar{y},0)$ for $\bar{y}\in S(\bar{x})$. Then
	\begin{equation}
	\emptyset \neq \limsup_{x\to\bar{x},\mu\to 0} \nabla \alert{\vmuquad} (x)\subseteq \bigcup _{\bar{y}\in S(\bar{x})} \partial_ x g(\bar{x},\bar{y}).
	\label{eq:almost_gradientconsistent}
	\end{equation} 
\end{lemma}

\ifdefined\submit
\noindent {\textbf{Proof.}}
\else
\begin{proof}
\fi
	Consider any sequence $\mu_k\searrow 0$ and $ x^k\to \bar{x}$. That $ \limsup_{x\to\bar{x},\mu\to 0} \nabla \alert{\vmuquad} (x)$ is non-empty holds due to our boundedness assumption. Suppose that $ \nabla v_{\mu_k}(x^k)\to \xi$. By \cref{prop:vmu_C1}, $\nabla_x g_{\mu_k}(x^k, y^k)\to \xi$ where $y^k= S_{\mu_k}(x^k)$. Using \eqref{eq:v_intermsof_h}, the uniform level boundedness of $h$ given by \eqref{eq:h_extended}, and \cref{thm:uniformlevelbounded}(b), we know that $\{y^k\}$ is bounded with accumulation points in $S(\bar{x})$, and so we may assume without loss of generality that $y^k\to \bar{y}\in S(\bar{x})$. 
	Meanwhile, using again our boundedness assumption, note that $\nabla g_{\mu_k}(x^k, y^k)$ is bounded so that we may assume without loss of generality that it converges to a point $(\xi,\eta)$, which belongs to $ \partial g(\bar{x},\bar{y})$ by gradient consistency of $\{\gmu :\mu>0 \}$. With \eqref{eq:projection_first_element_subset}, we infer that $\xi\in \partial_x g(\bar{x},\bar{y})$ and therefore \eqref{eq:almost_gradientconsistent} holds. 
\ifdefined\submit
\hfill \Halmos
\else 
\end{proof}
\fi 

In general, under some regularity conditions (see \cite[Definition 7.25]{RW98}) and the assumption that $g(x,\cdot)$ is convex for all $x\in O_X$, we can only obtain 	
\begin{equation}
\partial v(\bar{x}) \subseteq \bigcup_{\bar{y}\in S(\bar{x})} \partial_x g(\bar{x},\bar{y})
\label{eq:valuefunction_convexobjective2}	
\end{equation}
by using \cite[Theorem 10.13]{RW98}, \cite[Proposition 2.3.15]{Clarke83}, and \cite[Corollary 1 of Theorem 2.9.8]{Clarke83}. Hence, \eqref{eq:almost_gradientconsistent} is in general a weaker result as it only guarantees that $ \limsup_{x\to\bar{x},\mu\to 0} \nabla \alert{\vmuquad} (x)$ is contained on the right hand side of \eqref{eq:valuefunction_convexobjective2}. Nevertheless, with the above, we immediately obtain the gradient consistency of $\{\alert{\vmuquad} :\mu >0\}$ for several special cases.

\begin{theorem}
	\alert{Suppose that \cref{assume:quadreg,assume:quadred_2,assume:quadred_3,assume:quadred_4} hold.} and that $\{\gmu: \mu>0\}$ is a smoothing family of functions satisfying the gradient consistent property. Let $\bar{x}\in O_X$ and suppose that $\{\nabla \gmu (x,y) : (x,y,\mu) \in \Theta\}$ is bounded for some neighborhood $\Theta$ of $(\bar{x},\bar{y},0)$ for $\bar{y}\in S(\bar{x})$. Then, under any of the following conditions, $\{\alert{\vmuquad} : \mu>0\}$ is a family of smooth approximations of $v$ that satisfies the gradient consistent property at $\bar{x}$:
	\begin{enumerate}[(a)]
		\item $\nabla_x g(x,y)$ exists and is continuous on $\Omega\times Y$; 
		\item $g(\cdot, y)$ is $\rho$-weakly concave on $\Omega$ for all $y\in O_Y$; or
		\item \eqref{eq:partialformula} holds and $g$ is convex in $(x,y)$,
	\end{enumerate}
	where $\Omega$ is a neighborhood of $\bar{x}$. 
\end{theorem}
\ifdefined\submit
\noindent {\textbf{Proof.}}
\else
\begin{proof}
\fi
	We first show that condition \eqref{eq:projection_first_element_subset} holds for each of the cases described. To this end, let $\xi \in \pi_x (\partial g(\bar{x},\bar{y})) $ and let $\eta\in \Re^m$ such that $(\xi,\eta)\in \partial g(\bar{x},\bar{y})$. For (a), it is enough to show that $\xi = \nabla_x g(\bar{x},\bar{y})$. Let $\{ (x^k,y^k)\}$ be a sequence converging to $(\bar{x},\bar{y})$ such that $g$ is differentiable at each $(x^k,y^k)$ and $\{\nabla g(x^k,y^k)\}$ is convergent. By continuity of $\nabla_x g$, we conclude that $\nabla_x g(x^k,y^k) \to \nabla_x g(\bar{x},\bar{y})$. It now directly follows from \cite[Theorem 2.5.1]{Clarke83} that  $\xi = \nabla_x
	g(\bar{x},\bar{y})$, proving that \eqref{eq:projection_first_element_subset} holds. For (b), we have from \cite[Proposition 2.5.3]{Clarke83} that $\pi_x \left( \partial \hat{h}(\bar{x},\bar{y})\right) \subseteq \partial_x \hat{h}(\bar{x},\bar{y})$ where $\hat{h}(x,y) \coloneqq h(x,y) + \frac{\rho}{2}\|x\|^2$ and $h=-g$ as in the proof of \cref{thm:danskin_extended_noncompact}. The inclusion \eqref{eq:projection_first_element_subset} immediately follows. Lastly, in part (c), the hypothesis that \eqref{eq:partialformula} holds directly implies \eqref{eq:projection_first_element_subset}. Hence, we may use \cref{thm:gradientconsistentproperty}. By \cref{thm:danskin_noncompact}, \cref{thm:danskin_extended_noncompact} and 
	\cref{thm:danskin_convex}, which can be applied for cases (a), (b) and (c), respectively, the right-hand side of \eqref{eq:almost_gradientconsistent} is either equal to or contained in $\partial v(\bar{x})$. Hence, gradient consistency holds at $\bar{x}$.
\ifdefined\submit
\hfill \Halmos
\else 
\end{proof}
\fi

\section{Entropic Regularization}\label{sec:entropy} \cite{LinXuYe14} proposed the function $ x\mapsto -\mu \ln \left( \int_Y \exp \left( -\mu^{-1} g(x,y)\right) dy \right)$ as a smooth approximation\footnote{For an intuitive interpretation of this smoothing function, recall that a usual approximation of $\max \{ x_1,x_2,\dots, x_r\}$ is $\mu \ln \sum_{i=1}^r \exp (\mu^{-1}x_i) $. We may thus view the given smoothing function as a continuous extension of this approximation.} of the value function when $g$ is smooth and $Y$ is a \textit{compact} set. This is based on the earlier works of \cite{FangWu96,LiFang97} on smoothing by entropic regularization for the maximum value function. A natural way to extend this to a nonsmooth function $g$ is to replace $g$ with a smooth approximation. That is, if $\{ \gmu:\mu >0\}$ is a family of smooth approximations of $g$, then we consider
\begin{equation}
\alert{\vmuentrop} (x) \coloneqq -\mu \ln \left( \int_Y \exp \left( -\mu^{-1}\gmu(x,y)\right) dy \right)
\label{eq:smoothing_extended_linxuye}
\end{equation}
as a smooth approximation of the value function. 

We show that when $Y$ is a compact set, the function \eqref{eq:smoothing_extended_linxuye} is indeed a smooth approximation of $v$ which satisfies gradient consistency under certain conditions. In fact, we can further extend this result to the case when $Y$ is non-compact by replacing $Y$ in \eqref{eq:smoothing_extended_linxuye} with $\Ymu$,
%\[ \alert{\vmuentrop} (x) \coloneqq -\mu \ln \left( \int_{\Ymu} \exp \left( -\mu^{-1}\gmu(x,y)\right) dy \right),\]
where $\{\Ymu : \mu >0 \}$ is a family of compact subsets of $Y$ such that $\Ymu\nearrow Y$, \alert{that is, $
Y_\mu \subseteq Y_{\mu'} \subseteq Y$  for all $0 < \mu' < \mu$, and $\bigcup_{\mu > 0} Y_\mu = Y.$} These considerably extend the smoothing framework of \cite{LiFang97,LinXuYe14} to unbounded closed sets and nonsmooth functions. Such an extension is not trivial, neccessitating subtle modifications to the proofs established for the compact case, as well as sufficient conditions required to obtain a smoothing family and gradient consistency. These intricacies are systematically discussed in this manuscript. To concentrate on the main ideas, we focus first on the compact case in \cref{sec:entropy_compact}, and then present analogous results for the noncompact case in \cref{sec:entropy_noncompact}. 

\subsection{Entropic regularization for bounded constraint set}\label{sec:entropy_compact}

We start by proving that $\alert{\vmuentrop}$ given by \eqref{eq:smoothing_extended_linxuye} is continuously differentiable.  

\begin{proposition}
	\label{prop:smoothing_entropy}
	Suppose that $Y$ is compact and $\{\gmu : \mu >0\}$ is a family of smooth approximations of $g$ over some open set $O_X\times O_Y$. Then for each $\mu>0$, $\alert{\vmuentrop}$ given by \eqref{eq:smoothing_extended_linxuye} is continuously differentiable on $O_X$.
\end{proposition}
\ifdefined\submit
\noindent {\textbf{Proof.}}
\else
\begin{proof}
\fi 
	Let $\Omega$ be any bounded open set whose closure is contained in $O_X$. Note that by the continuity of $\gmu$, $\| \nabla_x \exp \left( -\mu^{-1}\gmu(x,y)\right)\| $ is uniformly bounded over  $\cl (\Omega)\times Y$ by any constant $c$ such that 
	\[ c\geq \max_{(x,y)\in \cl (\Omega)\times Y} \mu^{-1} \exp (\mu^{-1}\gmu(x,y)) \| \nabla_x \gmu (x,y)\|,\]	
	By Leibniz rule, we obtain
	\begin{equation}
	\nabla \alert{\vmuentrop} (x) = \int_Y \alpha_{\mu}(x,y) \nabla_x \gmu (x,y) dy \quad \forall x\in \Omega,
	\label{eq:vmu_integral}
	\end{equation}
	where 
	\begin{equation}
	\alpha_{\mu} (x,y) \coloneqq \frac{\exp(-\mu^{-1}\gmu (x,y)) }{\int_Y \exp(-\mu^{-1}\gmu (x,z))  dz}.
	\label{eq:alpha_mu}
	\end{equation}
	We now show that $\nabla \alert{\vmuentrop}$ is continuous on $\Omega$. First, we claim that $\alpha_{\mu}$ is a continuous function on $\Omega\times Y$. Note that $\exp(-\mu^{-1}\gmu (x,z))$ is uniformly bounded above by some constant over the compact set $\cl (\Omega)\times Y$. Hence, continuity of $\alpha_{\mu}$ on $\Omega\times Y$ holds by the bounded convergence theorem. By continuity of $\alpha_{\mu}$ and $\nabla_x \gmu$, we can again apply the same argument to show that $\nabla \alert{\vmuentrop}$ is continuous on $\Omega$. Since $\Omega$ is arbitrarily chosen, we have completed the proof. 
\ifdefined\submit
\hfill \Halmos
\else 
\end{proof}
\fi 

%\begin{remark}
%	Note that the compactness of $X$ may be replaced with any other condition that will guarantee the existence of a function $V_1,V_2,V_3\in L^1(Y)$ such that for any $x$, $\| \nabla_x \exp \left( -\mu^{-1}\gmu(x,\cdot)\right)\|$, $\exp(-\mu^{-1} \gmu (x,\cdot)) $, and $\alpha_{\mu}(x,\cdot) \| \nabla_x \gmu (x,\cdot)\| $ are bounded above by $V_1, V_2$ and $V_3$.
%\end{remark}

\alert{Similar to \cref{prop:ennamcq_quad}, we show that ENNAMCQ is satisfied at any point $(\bar{x},\bar{y})\in X\times Y$ for \eqref{eq:bilevel_nlp_formulation_epsilon_mu} with $v_{\mu} = \vmuentrop$.} 
\alert{\begin{proposition}\label{prop:ennamcq_entrop_compact}
		Suppose that \cref{assume:quadreg}(b) holds, and let $v_{\mu}=\vmuentrop$ in \eqref{eq:bilevel_nlp_formulation_epsilon_mu} with $\epsilon>0$. For any $\mu>0$ such that $\mu \ln \vol (Y) < \epsilon$, ENNAMCQ holds at any point $(\bar{x},\bar{y})\in X\times Y$. 
%		, that is, either $g_{\mu}(\bar{x},\bar{y})-\vmuentrop (\bar{x})<\epsilon$ or $g_{\mu}(\bar{x},\bar{y})-\vmuentrop (\bar{x})\geq \epsilon$ but $0\notin \nabla g_{\mu}(\bar{x},\bar{y}) - \nabla \vmuentrop(\bar{x})\times \{0\} + N_{X\times Y}(\bar{x},\bar{y}). $
\end{proposition}}
\alert{\ifdefined\submit
	\noindent {\textbf{Proof.}}
	\else
	\begin{proof}
		\fi
		Suppose $g_{\mu}(\bar{x},\bar{y})-\vmuentrop (\bar{x})\geq \epsilon$, and  assume, to the contrary, that $0\in \nabla g_{\mu}(\bar{x},\bar{y}) - \nabla \vmuentrop(\bar{x})\times \{0\} + N_{X\times Y}(\bar{x},\bar{y})$. As in the proof of \cref{prop:ennamcq_quad}, we have $\bar{y}\in \argmin_{y\in Y} g_{\mu}(\bar{x},y)$. Thus, $g_{\mu}(\bar{x},\bar{y}) \leq g_{\mu}(\bar{x},y)$ for all $y\in Y$, which consequently implies that 
		\[-\mu \ln \left( \int_Y \exp (-\mu^{-1} g_{\mu}(\bar{x},\bar{y}))dy  \right) \leq \vmuentrop (\bar{x}) \leq g_{\mu}(\bar{x},\bar{y})-\epsilon .\]
		It can be shown that this implies that $\mu \ln \vol (Y) \geq \epsilon$, which contradicts our choice of $\mu$. Therefore, $0\notin \nabla g_{\mu}(\bar{x},\bar{y}) - \nabla \vmuentrop(\bar{x})\times \{0\} + N_{X\times Y}(\bar{x},\bar{y})$, as desired. 
		\ifdefined\submit
		\hfill \Halmos
		\else 
	\end{proof}
	\fi }

We now show that $\{ \alert{\vmuentrop}: \mu >0\}$ is a family of smooth approximations in the sense of \cref{defn:smoothing}. We need the following lemma.% whose proof is given in \cref{appendix:proof_integral_lowerbound}.

\begin{lemma}
	\label{lemma:integral_lower_bound}
	Let $Y$ be a compact set and let $\Omega$ be any convex open set such that $\cl (\Omega) \subseteq O_X$.  For any $\tau \in (0,1)$, there exists $\delta>0$ such that for any $\mu\in (0,\delta)$, 
	\begin{equation}
	\tau (\mu \vol(Y))^{\mu} \max _{y\in Y} \exp (-g(x,y)) \leq \left( \int_Y \exp \left( -\mu^{-1}g(x,y)\right)dy \right) ^{\mu} \quad \forall x\in \cl (\Omega).
	\label{eq:integral_lower_bound}
	\end{equation}
\end{lemma}
\ifdefined\submit
\noindent {\textbf{Proof.}}
\else
\begin{proof}
\fi 
	From the proof of \cite[Theorem 1]{FangWu96}, the result follows if we can show that for any $\tau \in (0,1)$, there exists a sufficiently large $N$ such that 
	\begin{equation}
	\ds \max_{j=1,2,\dots,k} ~~\min_{y\in Y_j} \exp (-g(x,y)) > \tau \cdot \max_{y\in Y} \exp (-g(x,y)) \quad \forall x\in \cl (\Omega)
	\label{eq:maxminmaxinequality}
	\end{equation}
	holds for all integers $k\geq N$ and any partition $\{Y_1,Y_2,\dots,Y_k\}$ of $Y$ with $\max_{j=1,2,\dots,k} \diam (Y_j)\to 0$ as $k\to \infty$, where $\diam (Y_j)$ is the diameter of $Y_j$. This is similar to \cite[Lemma 2]{FangWu96} where the inequality \eqref{eq:maxminmaxinequality} holds on the set $X$ but the function $g$ should be such that $g(\cdot,y)$ is continuously differentiable for any $y\in Y$. We show that in fact, the said result can be extended to our setting that $g$ is merely Lipschitz continuous.
	
	For the sake of contradiction, suppose that there exists $\tau\in (0,1)$, a sequence $\{x^k\}\subseteq \cl (\Omega)$, and integers $i_k$ with $i_k\to \infty$ such that 
	\[ \frac{\ds \max_{j=1,2,\dots,i_k} ~~\min_{y\in Y_j} \exp (-g(x,y))}{\tau \cdot \max_{y\in Y} \exp (-g(x,y))}\leq 1\]
	for some partition $\{Y_1,Y_2,\dots,Y_{i_k} \}$ with $\max_{j=1,2,\dots,i_k} \diam (Y_j)\to \infty$ as $k\to \infty$. By the compactness of $\cl (\Omega)$, we may assume that $x^k\to x^*$ for some $x^*\in \cl (\Omega)$. Now, we use the mean-value theorem for Lipschitz continuous functions. Noting that $\cl (\Omega )\subseteq O_X$, we have from \cite[Theorem 2.3.7]{Clarke83}  that there exists $u^k$ in the relative interior of the line segment joining $x^*$ and $x^k$ such that 
	\[ \exp (-g(x^k,y)) - \exp (-g(x^*,y)) \in \lla \partial h_y(u^k),x^k-x^*\rla ,\]
	where $h_y(x)\coloneqq \exp(-g(x,y))$. By Theorem 2.3.9(ii) and Proposition 2.3.1 of \cite{Clarke83}, we have
	\begin{equation}
	\exp (-g(x^k,y)) - \exp (-g(x^*,y)) = -\lla \exp (-g(u^k,y))v^k,x^k-x^*\rla ,
	\label{eq:taylor}
	\end{equation}
	where $v^k\in \partial _xg(u^k, y)$. Meanwhile, noting the convexity of $\cl (\Omega)$, each $u^k$ is in $\cl (\Omega)\subseteq O_X$. On the other hand, $\partial g(u,y)$ is a compact set for any $(u,y)$ by \cite[Proposition 7.1.4(a)]{FP03}. Hence, with the upper-semicontinuity of the Clarke subdifferential (see \cite[Proposition 7.1.4(b)]{FP03} or \cite[Proposition 2.1.5(d)]{Clarke83}), we conclude that $\bigcup_{(u,y)\in \cl (\Omega)\times Y} \partial g(u,y)$ is bounded. (Note that here, the Lipschitz continuity of $g$ over the open set $O_X\times O_Y$ containing the compact set $\cl (\Omega)\times Y$ plays an important role; see also \cref{sec:prelim_nonsmooth}). Since $\partial_x g(u,y) \subseteq \pi_x \left( \partial g (u,y)\right)$ by \cite[Proposition 2.3.16]{Clarke83}, we conclude that the set $\bigcup_{(u,y)\in \cl (\Omega)\times Y} \partial_x g(u,y)$ is likewise bounded, which together with the continuity of $g$ implies that 	$M\coloneqq \sup \left\lbrace \exp (-g(u,y))\|v\| :v\in  \partial_x g(u,y)~\text{and}~ (u,y)\in\cl(\Omega)\times Y\right\rbrace$ is finite\jhsolved{Yes, the ``$M$'' is the same for all $k$. I added the definition of $M$ to make this clearer.}. Going back to \eqref{eq:taylor} and using Cauchy-Schwarz inequality, we have
	\[ \min_{y\in Y_j} \alert{\exp }(-g(x^*,y))  - M \|x^k-x^*\| \leq \min_{y\in Y_j} \exp (-g(x^k,y)) \leq  \min_{y\in Y_j}\exp (-g(x^*,y)) + M\|x^k-x^*\|. \]
	At this point, we follow arguments analogous to those in the proof of \cite[Lemma 2]{FangWu96} to obtain a contradiction\jhsolved{The contradiction you were looking for before is obtained by following the same arguments as in the reference indicated.}, thus completing the proof. We note, however, that showing that equation (17) in \cite{FangWu96} does not require the ``super uniform continuity'' assumption imposed by the authors, as it automatically follows from \cref{thm:uniformlevelbounded}. 
	%	In conclusion, we have shown that for any $\tau\in (0,1)$ and any subset $\Omega$ with $\cl (\Omega) \subseteq O_X$, there exists a sufficiently large $N=N(\Omega)$ such that  
	%		\begin{equation}
	%	\frac{\ds \max_{j=1,2,\dots,k} ~~\min_{y\in Y_j} \exp (-g(x,y))}{\tau \cdot \max_{y\in Y} \exp (-g(x,y))}>1 \quad \forall x\in \cl (\Omega)
	%	\label{eq:maxminmaxinequality2}
	%	\end{equation}
	%	holds for all $k\geq N$ and any partition $\{Y_1,Y_2,\dots,Y_k\}$ of $Y$ with $\max_{j=1,2,\dots,k} \diam (Y_j)\to 0$ as $k\to \infty$, where $\diam (Y_j)$ is the diameter of $Y_j$. Since $O_X$ is bounded set, it can be covered
\ifdefined\submit
\hfill \Halmos
\else 
\end{proof}
\fi

\begin{theorem}
	\label{thm:smoothing_integral}
	\alert{Suppose that \cref{assume:quadred_4} holds}. Then $\{\alert{\vmuentrop} : \mu >0\}$ is a family of smooth approximations over $O_X$ for the value function $v$.
\end{theorem}

\ifdefined\submit
\noindent {\textbf{Proof.}}
\else
\begin{proof}
\fi 
	\cref{prop:smoothing_entropy} already provides the continuous differentiability of each $\alert{\vmuentrop}$ on $O_X$. We just need to show that for any $\bar{x}\in O_X$, $\alert{\vmuentrop}(x)\to v(\bar{x})$ as $\mu\to 0$ and $x\to \bar{x}$. Take any $\bar{x}\in O_X$ and let $\varepsilon >0$ be given, and $\tau \in (\exp (-\varepsilon),1)$. By \cref{lemma:integral_lower_bound}, given a convex open set $\Omega$ containing $\bar{x}$ such that $\cl(\Omega) \subseteq O_X$, there exists $\delta>0$ such that \eqref{eq:integral_lower_bound} holds for all $\mu\in (0,\delta)$. Together with \cite[Theorem 1]{FangWu96}, we have
	\[\tau (\mu \vol(Y))^{\mu} \max _{y\in Y} \exp (-g(x,y)) \overset{\eqref{eq:integral_lower_bound}}{\leq} \left( \int_Y \exp \left( -\mu^{-1}g(x,y)\right) dy\right) ^{\mu}\leq \vol(Y)^{\mu} \max_{y\in Y} \exp (-g(x,y)) \]
	for all $x\in \cl (\Omega)$. Noting the monotonicity of the logarithmic function, we have
	\begin{multline*}
	- \mu \ln \vol(Y) -  \max_{y\in Y} (-g(x,y)) \leq -\mu \ln \left( \int_Y \exp \left( -\mu^{-1}g(x,y)\right) dy\right) \\ \leq -\ln \tau - \mu \ln  (\mu \vol(Y))-  \max _{y\in Y} (-g(x,y)).
	\end{multline*}
	By the choice of $\tau$ and the definition of $v$, we have
	\begin{equation*}
	v(x) - \mu \ln \vol(Y) \leq -\mu \ln \left( \int_Y \exp \left( -\mu^{-1}g(x,y)\right) dy\right)  \leq v(x) - \mu \ln  (\mu \vol(Y)) + \varepsilon.
	\end{equation*}
	We may choose $\delta>0$ to be sufficiently small enough so that for all $\mu\in (0,\delta)$, we have $ \mu \ln \vol(Y)<2\varepsilon$ and $ - \mu \ln  (\mu \vol(Y))< \varepsilon $. Hence, 
	\begin{equation}
	\left| -\mu \ln \left( \int_Y \exp \left( -\mu^{-1}g(x,y)\right) dy\right) - v(x) \right| < 2\varepsilon,
	\label{eq:bounds_v(x)}
	\end{equation}
	valid for all $x\in \cl(\Omega)$ and $\mu\in (0,\delta)$. On the other hand, invoking the uniform convergence of $\gmu$ to $g$ on $\cl (\Omega)\times Y$, we may choose $\delta$ to be small enough so that \eqref{eq:bounds_v(x)} holds and $|\gmu (x,y) - g(x,y)|<\varepsilon$ for all $\mu \in (0,\delta)$ and $(x,y)\in \cl (\Omega)\times Y$. \alert{Consequently, $\gmu (x,y) - \varepsilon < g(x,y) < \gmu(x,y)+\varepsilon$, and therefore 
	\[\exp (-\mu^{-1} \gmu(x,y)) \exp (\mu^{-1}\varepsilon) > \exp (-\mu^{-1}g(x,y)) > \exp (-\mu^{-1}\gmu(x,y)) \exp (-\mu^{-1}\varepsilon), \]
	for all $x\in \cl (\Omega)$ and $\mu \in (0,\delta)$. Integrating over $Y$ and noting the monotonicity of logarithmic function, we obtain }
	\begin{equation*}
	\left| -\mu \ln \left( \int_Y \exp \left( -\mu^{-1}\gmu (x,y)\right) dy\right) + \mu \ln \left( \int_Y \exp \left( -\mu^{-1}g(x,y)\right) dy\right) \right| <  \varepsilon
	\end{equation*}
	for all $x\in \cl (\Omega)$ and $\mu \in (0,\delta)$. Together with \eqref{eq:bounds_v(x)} and the definition of $\alert{\vmuentrop}$, we have 
	\begin{equation}
	| \alert{\vmuentrop} (x) - v(x) | < 3\varepsilon \qquad \forall x\in \cl (\Omega), ~\forall \mu \in (0,\delta) .
	\label{eq:bounds_vmu}
	\end{equation}
	By the continuity of the value function on $ \cl (\Omega)$ (see \cref{thm:uniformlevelbounded}(a)), we may choose $\delta>0$ so that $|v(x) - v(\bar{x})|< \varepsilon$ for all $x\in \cl (\Omega)$ with $\|x - \bar{x}\|<\delta$. Together with \eqref{eq:bounds_vmu}, we have $|\alert{\vmuentrop} (x) - v(\bar{x}) |< 4\varepsilon$ for all $\mu\in (0,\delta)$ and $x\in \cl (\Omega)$ with $\|x-\bar{x}\|<\delta$. This completes the proof. 
\ifdefined\submit
\hfill \Halmos
\else 
\end{proof}
\fi 

%\begin{remark}
%	There are several of smoothing family $\{ \gmu \}$ that converges uniformly to $g$ on their given domain. For instance, when $g(\cdot,y)$ is smooth for any $y\in O_Y$ and $g(x,\cdot)$ is convex for any $x\in O_X$, then the smoothing functions via Moreau envelope satisfies the uniform convergence property; see \cite[Theorem 10.51]{Beck17}. 
%\end{remark}

We next show that $\alert{\vmuentrop}$ satisfies the gradient consistent property. To this end, we need the following lemmas. Note that \cref{lemma:alpha_outside_solutionset} is an extension of \cite[Theorem 5.3]{LinXuYe14}, but we provide a simpler and shorter proof as shown below. Moreover, we also need another technical result shown in \cref{lemma:alpha_outside_solutionset_Integral}, which will be helpful in our proof of gradient consistency. 

\begin{lemma}
	\label{lemma:alpha_outside_solutionset}
	\alert{Suppose that \cref{assume:quadred_4} holds. If $\bar{x}\in O_X$ and $y'\in Y\setminus S(\bar{x})$, then } $\alpha_{\mu}(x,y')\to 0$ as $(x,\mu)\to (\bar{x},0)$, where $\alpha_{\mu}$ is given by \eqref{eq:alpha_mu}.
\end{lemma}

\ifdefined\submit
\noindent {\textbf{Proof.}}
\else
\begin{proof}
\fi 
	By \cref{thm:smoothing_integral}, $-\gmu(x,y')+\alert{\vmuentrop} (x) \to -g(\bar{x},y')+v(\bar{x})$ as $(x,\mu)\to (\bar{x},0)$. Meanwhile, 	observe that we can rewrite $\alpha_{\mu}$ as $\alpha_{\mu} (x,y) = \exp \left( \mu^{-1} (-\gmu (x,y)+\alert{\vmuentrop} (x))\right)  .$
	Hence, letting $(x,\mu)\to (\bar{x},0)$ and noting that $-g(\bar{x},y')+v(\bar{x})<0$ since  $y'\notin S(\bar{x})$, we see that $\alpha_{\mu}(x,y')\to 0$, as desired.  
\ifdefined\submit
\hfill \Halmos
\else 
\end{proof}
\fi 

\begin{lemma}
	Let $\bar{x}\in O_X$ and let $Z$ be any closed subset of  $Y\setminus S(\bar{x})$. \alert{If \cref{assume:quadred_4} holds, then} $\int_{Z} \alpha_{\mu}(x,y) dy\to 0$ as $(x,\mu) \to (\bar{x},0)$.
	\label{lemma:alpha_outside_solutionset_Integral}
\end{lemma}
\ifdefined\submit
\noindent {\textbf{Proof.}}
\else
\begin{proof}
\fi 
	First, we claim that $\gmu(x,\cdot) \to g(\bar{x},\cdot)$ as $(x,\mu)\to (\bar{x},0)$ uniformly on $Y$. Note that by hypothesis, $\gmu\to g$ as $\mu\to 0$ uniformly on $\cl (B_{\delta_0}(\bar{x}))\times Y$, where $\delta_0$ is any positive number such that $\cl (B_{\delta_0}\alert{(\bar{x})})\subseteq O_X$. Thus, given any $\varepsilon>0$, there exists $\delta \in (0,\delta_0)$ such that 
	\begin{equation}
	\mu\in (0,\delta) \quad \Longrightarrow \quad |\gmu (x,y)-g(x,y)|<\varepsilon/2\quad \forall  (x,y) \in \cl (B_{\delta_0}(\bar{x}))\times Y.
	\label{eq:uniconvergent}
	\end{equation}
	Meanwhile, since $g$ is uniformly continuous on the compact set $\cl (B_{\delta_0}(\bar{x}))\times Y$, we may choose a smaller $\delta>0$, if necessary, so that 
	\begin{equation}
	|g(x,y) - g(x',y')|<\varepsilon/2 \quad \forall (x,y),(x',y')~\text{with}~\norm{(x,y)-(x',y')}<\delta.
	\label{eq:unicontinuous}
	\end{equation}
	Since $\delta<\delta_0$ and $\norm{(x,y)-(\bar{x},y)}=\norm{x-\bar{x}}$, then for any $(x,\mu)$ such that $\norm{x-\bar{x}}<\delta$ and $\mu \in (0,\delta)$, it follows from triangle inequality and \eqref{eq:uniconvergent}-\eqref{eq:unicontinuous} that $|\gmu(x,y) -g(\bar{x},y)|<\varepsilon$ for any $y\in Y$. This completes that $\gmu(x,\cdot) \to g(\bar{x},\cdot)$ as $(x,\mu)\to (\bar{x},0)$ uniformly on $Y$.

	Now, let $\varepsilon \coloneqq \ds \max_{y\in Z} g(\bar{x},y)-v(\bar{x})$. Note that since $Z$ is compact and lies outside $S(\bar{x})$, then $\varepsilon$ is a finite positive number. Since $\lim_{(x,\mu)\to (\bar{x},0)}\alert{\vmuentrop} (x) = v(\bar{x})$ by \cref{thm:smoothing_integral}, there exists $\delta>0$ such that $\alert{\vmuentrop}(x) - v(\bar{x})<\varepsilon/2$ for all $(x,\mu) \in B_{\delta}(\bar{x},0)$. On the other hand, by the claim proved above, we may choose a smaller $\delta>0$ if necessary so that $g(\bar{x},y) - \gmu(x,y) < \varepsilon/2$ for all $y\in Y$ and for all $(x,\mu) \in B_{\delta}(\bar{x},0)$. Consequently, $-\gmu(x,y) + \alert{\vmuentrop}(x) < -g(\bar{x},y)+v(\bar{x})+ \varepsilon$ for all $y\in Y$ and for all $(x,\mu) \in B_{\delta}(\bar{x},0)$. On the other hand, by our choice of $\varepsilon$, we have  $-g(\bar{x},y)+v(\bar{x})+ \varepsilon\leq 0$ for all $y\in Z$, so that $-\gmu(x,y) + \alert{\vmuentrop}(x) <0$ for all $(x,\mu)\in B_{\delta}(\bar{x},0)$ and $y\in Z$. In turn, $\alpha_{\mu}(x,\cdot)=\exp \left( \mu^{-1} (-\gmu (x,\cdot)+\alert{\vmuentrop} (x))\right)  $ is uniformly bounded above by 1 on $Z$ for all $(x,\mu)\in B_{\delta}(\bar{x},0)$. Finally, we obtain the desired limit $\int_{Z} \alpha_{\mu}(x,y) dy\to 0$ as a direct consequence of \cref{lemma:alpha_outside_solutionset} and the bounded convergence theorem.
\ifdefined\submit
\hfill \Halmos
\else 
\end{proof}
\fi 

We now prove the following result which is instrumental in showing gradient consistency, later proved in \cref{cor:gradconsistency_integral}. Note that unlike the smooth analogue \cite[Theorem 5.5]{LinXuYe14}, we do not assume compactness of the set $X$ in \cref{thm:gradconsistency_integral}. As mentioned above, the proof of the following lemma also needs the additional result \cref{lemma:alpha_outside_solutionset_Integral}. \alert{We also emphasize that, due to the nonsmoothness of $g$, the main challenge lies in finding an element of the set
	\alert{\begin{equation}\label{eq:P(x)}
	P(\bar{x}) \coloneqq \co \{\xi \in \Re^n: \xi \in \partial_x g(\bar{x},\bar{y})~\text{and}~\bar{y}\in S(\bar{x}) \}.
	\end{equation}}that is arbitrarily close to the set $\limsup_{x\to \bar{x}, \mu\searrow 0} \nabla \alert{\vmuentrop} (x)$}. Moreover, the subsequent verification that the identified element is indeed close requires some subtle techniques. \alert{The following assumption will also be needed. }
\alert{\begin{assumption}
	\label{assume:strongergradconsistency}
	At a given point $\bar{x}\in O_X$, $\dist (\nabla_x \gmu (x,\cdot), \partial_x g(x,\cdot)) $ converges to $0$ uniformly on $Y$ as $(x,\mu)\to (\bar{x},0)$.
\end{assumption}}

\begin{lemma} \alert{Let $\bar{x}\in O_X$. Suppose that \cref{assume:quadred_4,assume:strongergradconsistency} hold}. If there exists a neighborhood $\Omega$ of $\bar{x} $ such that $\partial_x g(\cdot,\cdot)$ is upper semicontinuous on $\Omega\times O_Y$, then 
	\[\emptyset \neq \limsup_{x\to \bar{x}, \mu\searrow 0}\nabla \alert{\vmuentrop} (x) \subseteq P(\bar{x}) ,\]
where $P(\bar{x})$ is given by \alert{\eqref{eq:P(x)}}. 
\label{thm:gradconsistency_integral}
\end{lemma}
\ifdefined\submit
\noindent {\textbf{Proof.}}
\else
\begin{proof}
\fi 
	First, we show that $\limsup_{x\to \bar{x}, \mu\searrow 0}\nabla \alert{\vmuentrop} (x)$ is indeed nonempty. By \alert{\cref{assume:strongergradconsistency}}, there exists $\delta'>0$ such that $\dist ( \nabla_x \gmu (x,y), \partial_x g(x,y)) < 1 $ for all $\mu \in (0,\delta') $ and $(x,y)\in \cl (B_{\delta'}(\bar{x}))\times Y.$
	It follows that 
	\begin{equation}
	\|\nabla_x \gmu(x,y)\| \leq 1 + M \qquad \forall \mu \in (0,\delta') , ~\forall (x,y)\in B_{\delta'}(\bar{x})\times  Y
	\label{eq:nabla_gmu_bound}
	\end{equation}
	where 
	\begin{equation}
	M\coloneqq \sup \{ \|\xi\| : \xi\in \partial_x g(x,y)~\text{and}~ (x,y)\in \cl (B_{\delta'}(\bar{x}))\times Y\}.
	\label{eq:M}
	\end{equation}
	We may choose $\delta'>0$ to be sufficiently small so that $\cl (B_{\delta'}(\bar{x}))\subseteq \Omega$. Thus, $M<+\infty$ by the upper semicontinuity of $\partial_x g (\cdot,\cdot)$ and compactness of $\cl(B_{\delta'}(\bar{x}))\times Y$. From \eqref{eq:vmu_integral} and noting that $\int_Y \alpha_{\mu} (x,y)dy = 1$, we see that $\{ \nabla \alert{\vmuentrop} (x) : x\in B_{\delta'}(\bar{x}) , \mu\in (0,\delta')\}$ is a bounded set, thus proving the first claim. 
	
	To prove the inclusion, suppose that $\nabla v_{\mu_k}(x^k)\to \xi\in\Re^n$ for some sequences $x^k\to \bar{x}$ and $\mu_k\searrow 0$ as $k\to \infty$. We need to show that $\xi \in P(\bar{x})$. To this end, let $\varepsilon>0$ be given. We will find an element $\xi^k$ of $P(\bar{x})$ such that $\| \nabla v_{\mu_{k}}(x^{k}) - \xi^k\| <\varepsilon$ for sufficiently large $k$. For any $y\in S(\bar{x})$, we have from upper semicontinuity that there exists $\delta_y>0$ such that 
	\begin{equation}
	\partial_x g(x',y') \subseteq \partial_x g(\bar{x},y) + \frac{\varepsilon}{3} B_1(0)\quad \forall x'\in B_{\delta_y}(\bar{x}), ~y'\in B_{\delta_y}(y).
	\label{eq:usc}
	\end{equation}
	Meanwhile, since $Y$ is compact and $g$ is continuous, $S(\bar{x})$ is compact. Hence, there exist points $y_1,y_2,\dots, y_r\in S(\bar{x})$ and numbers $\delta_1,\delta_2,\dots,\delta_r>0$ such that \alert{$\delta_i=\delta_{y_i}$ and} $S(\bar{x})\subseteq \bigcup_{i=1}^r B_{\delta_i}(y_i)$. Denote $Y_1 \coloneqq Y\cap B_{\delta_1}(y_1)$, $Y_i  \coloneqq  Y \cap \left(B_{\delta_i} (y_i) \setminus \bigcup_{j=1}^{i-1}B_{\delta_j}(y_j) \right)$ for $i=2,\dots, r$, and $Y_{r+1} \coloneqq Y\setminus \bigcup_{i=1}^r B_{\delta}(y_i)$. 	Note that since $Y_{r+1} \subseteq Y \setminus S(\bar{x})$, it follows from \cref{lemma:alpha_outside_solutionset} that $\int_{Y_{r+1}} \alpha_{\mu_k}(x^k,y)  dy\to 0$. 
	
	\alert{Using again \cref{assume:strongergradconsistency},} we can find $\delta>0$ such that $\delta<\min \{ \delta', \delta_1,\dots, \delta_r\}$ and 
	\begin{equation}
	\dist ( \nabla_x \gmu (x,y), \partial_x g(x,y)) < \frac{\varepsilon}{3} \quad \forall \mu \in (0,\delta)  ~\text{and}~\forall (x,y)\in \cl (B_{\delta}(\bar{x}))\times Y.
	\label{eq:uniconvergence_derivative}
	\end{equation}
	By \cref{lemma:alpha_outside_solutionset_Integral}, we may choose $k$ large enough so that $x^k\in B_{\delta}(\bar{x})$, $\mu_k \in (0,\delta)$ and 
	\begin{equation}
	\int_{Y_{r+1}} \alpha_{\mu_k}(x^k,y)  dy < \frac{\varepsilon}{3(2M+1)} 
	\label{eq:integral_alpha_epsilon}
	\end{equation}
	where $M$ is given by \eqref{eq:M}. We now define $\xi^k$ as 
	\begin{equation*}
	\xi^k \coloneqq \sum_{i=1}^{r+1}\int_{Y_i} \alpha_{\mu_k}(x^k,y) \Pi_{\partial_xg(\bar{x},y_i)} \left( \varphi^{(i)}_k(y)\right)dy,
	\label{eq:xi^k}
	\end{equation*}
	where $y_{r+1}$ is an arbitrarily chosen element of $S(\bar{x})$, $\Pi_{\partial_xg(\bar{x},y_i)}$ is the projection operator onto the closed convex set $\partial_xg(\bar{x},y_i)$, and \alert{
	\begin{equation}
	\varphi^{(i)}_k (y)\coloneqq \Pi_{\partial _x g(x^k,y)} \left(\nabla_x g_{\mu_k}(x^k,y) \right)\quad \forall y\in Y_i.
	\label{eq:varphi_i}
	\end{equation}
	}Note that since $\bigcup_{i=1}^{r+1} Y_i = Y$, $\int_Y \alpha_{\mu_k} (x^k,y)dy = 1$, and the range of $\Pi_{\partial_xg(\bar{x},y_i)}$ is contained in $P(\bar{x})$, it follows that $\xi^k$ is indeed in $P(\bar{x})$ by Jensen's inequality. 
	
		%	By Rademacher's theorem, note that the Lebesgue measure of $E_i$ is zero since $g$ is Lipschitz continuous on $O_X\times O_Y$. Moreover, $\partial_x g(x^k,y) = \{ \nabla_x g(x^k,y)\}$ and $\varphi^{(i)}_k(y) = \nabla_x g(x^k,y)$ for any $y\in Y_i\setminus E_i$ where $i\in \{1,\dots,r\}$. Hence, for all $y\in Y_i\setminus E_i$ with $i\in \{1,\dots, r\}$,
	%
	
	We now show that $\| \nabla v_{\mu_{k}}(x^{k}) - \xi^k\| <\varepsilon$. 
%	For each $i=1,\dots, r$, denote $E_i \coloneqq \{y\in Y_i : g~\text{is not differentiable at }(x^k,y) \} $. 
\alert{Let $i\in \{1,\dots,r\}$ and $y\in Y_i$.	
First, we have from  \eqref{eq:uniconvergence_derivative} that 
	\begin{align}
	\norm{ \nabla_x g_{\mu_k}(x^k,y) - 	\varphi^{(i)}_k (y)} & = \norm{\nabla_x g_{\mu_k}(x^k,y) - \Pi_{\partial_x g(x^k,y)} (\nabla_x g_{\mu_k}(x^k,y))} \notag \\
	& = \dist(\nabla_x g_{\mu_k}(x^k,y),\partial_x g(x^k,y)) \notag \\
	& < \frac{\varepsilon}{3}.\label{eq:dist}
	\end{align}
	Meanwhile, note that $\varphi^{(i)}_k (y)\in \partial_x g(x^k,y)$. By \eqref{eq:usc} and noting that $y\in B_{\delta_i}(y_i)$ due to our choice of $\delta$, there exists $\nu_i(y) \in \partial_x g(\bar{x},y_i)$ such that $\norm{ \varphi^{(i)}_k (y) - \nu_i (y)} < \frac{\varepsilon}{3}$. It follows that $\norm{ \varphi^{(i)}_k (y) - \Pi_{\partial_xg(\bar{x},y_i)} \left( \nabla_x g(x^k,y)\right)} \leq \norm{ \varphi^{(i)}_k (y) - \nu_i (y)} <\frac{\varepsilon}{3}$. This, together with \eqref{eq:dist}, implies that 
	\begin{equation}
	\norm{ \nabla_x g_{\mu_k}(x^k,y) - \Pi_{\partial_xg(\bar{x},y_i)} \left( \varphi^{(i)}_k(y)\right)} < \frac{2\varepsilon}{3}, \quad \forall y\in Y_i, ~i\in \{1,\dots,r\},\label{eq:dist_S(x)_integral}
	\end{equation}
	by triangle inequality. Meanwhile, we also have
	\begin{align}
	&\norm{\int_{Y_{r+1}} \alpha_{\mu_k}(x^k,y) \left(\nabla_x g_{\mu_k}(x^k,y)- \Pi_{\partial_xg(\bar{x},y_{r+1})} \left( \varphi^{(r+1)}_k(y)\right)\right)dy } \notag \\
	&\quad \leq  \int_{Y_{r+1}} \alpha_{\mu_k}(x^k,y) \left( \norm{\nabla_x g_{\mu_k}(x^k,y)} + \norm{\Pi_{\partial_xg(\bar{x},y_{r+1})} \left( \varphi^{(r+1)}_k(y)\right)} \right)dy \notag \\
	& \quad \leq (2M+1) \int_{Y_{r+1}} \alpha_{\mu_k}(x^k,y)  dy \notag \\
	& \quad \overset{\eqref{eq:integral_alpha_epsilon}}{<} \frac{\varepsilon}{3}, \label{eq:outside_S(x)_integral}
	\end{align}
	where the second inequality holds by \eqref{eq:nabla_gmu_bound} and \eqref{eq:M} since $\delta<\delta'$. Finally, we have
	\begin{align}
	\| \nabla v_{\mu_k}(x^k) - \xi^k\| & = \norm{ \int_Y \alpha_{\mu_k}(x^k,y) \nabla_x g_{\mu_k} (x^k,y) dy - \sum_{i=1}^{r+1}\int_{Y_i} \alpha_{\mu_k}(x^k,y) \Pi_{\partial_xg(\bar{x},y_i)} \left( \varphi^{(i)}_k(y)\right)dy } \notag \\ 
	& \leq  \norm{  \sum_{i=1}^{r}\int_{Y_i} \alpha_{\mu_k}(x^k,y)\left( \nabla_x g_{\mu_k} (x^k,y)  -  \Pi_{\partial_xg(\bar{x},y_i)} \left( \varphi^{(i)}_k(y)\right) \right)dy} \notag \\
	& \qquad + \norm{\int_{Y_{r+1}} \alpha_{\mu_k}(x^k,y) \left(\nabla_x g_{\mu_k}(x^k,y)- \Pi_{\partial_xg(\bar{x},y_{r+1})} \left( \varphi^{(r+1)}_k(y)\right)\right)dy }\notag  \\
	& \overset{\eqref{eq:outside_S(x)_integral}}{<} \sum_{i=1}^{r}\int_{Y_i} \alpha_{\mu_k}(x^k,y)\norm{ \nabla_x g_{\mu_k} (x^k,y)  -  \Pi_{\partial_xg(\bar{x},y_i)} \left( \varphi^{(i)}_k(y)\right) }dy + \frac{\varepsilon}{3}\notag  \\
%	& = \sum_{i=1}^{r}\int_{Y_i\setminus E_i} \alpha_{\mu_k}(x^k,y)\norm{ \nabla_x g_{\mu_k} (x^k,y)  -  \Pi_{\partial_xg(\bar{x},y_i)} \left( \varphi^{(i)}_k(y)\right) }dy + \frac{\varepsilon}{3} \notag \\
	& \overset{\eqref{eq:dist_S(x)_integral}}{<} \frac{2\varepsilon}{3} \sum_{i=1}^{r}\int_{Y_i} \alpha_{\mu_k}(x^k,y) dy + \frac{\varepsilon}{3}\notag  \\
	& \leq \varepsilon . \label{eq:xi_epsilon_away}
	\end{align}}Since $\varepsilon>0$ is arbitrary, it follows that we can extract a subsequence $\{\nabla v_{\mu_k}(x^k) \}_{k\in K}$ where $K\subseteq \mathbb{N}$ together with $\{\xi^k\}_{k\in K}\subseteq P(\bar{x})$ such that $\|\nabla v_{\mu_k}(x^k) - \xi^k\| \to 0$ as $k\to \infty$ with $k\in K$. Since $P(\bar{x})$ is closed (due to the compactness of $S(\bar{x})$), any accumulation point $\xi^*$ of $\{\xi^k\}_{k\in K}$ belongs to $P(\bar{x})$. Since $\|\nabla v_{\mu_k}(x^k) - \xi^*\| \leq \|\nabla v_{\mu_k}(x^k) -\xi^k\| + \| \xi^k - \xi^*\|$ and $\nabla v_{\mu_k}(x^k) \to \xi$, we conclude that $\xi = \xi^* $, so that $\xi \in P(\bar{x})$ as desired. This completes the proof. 
\ifdefined\submit
\hfill \Halmos
\else 
\end{proof}
\fi 

\begin{remark}
	We note that the uniform convergence \alert{in \cref{assume:strongergradconsistency} which is required in the above lemma} can be achieved in many situations. Suppose that $g$ takes the form $g(x,y) = a(x,y) + b(p(x),q(y))$ where $a$ and $b$ are smooth, $b$ has Lipschitz continuous derivative, $p$ and $q$ are Lipschitz continuous functions, and the chain rule applies for $b(p(\cdot),q(y))$ for each $y$; for instance, under the conditions stipulated in \cite[Theorem 2.3.9]{Clarke83} or simply when $p$ is a smooth function. Let $\{ p_{\mu}: \mu>0\}$ and $\{ q_{\mu}: \mu>0\}$ be smooth approximations of $p$ and $q$, respectively, such that $\dist(\nabla p_{\mu}(\cdot),\partial p(\cdot)) \to 0$ pointwise (uniform convergence is not required), and $q_{\mu}\to q$ uniformly on $Y$. Note that these assumptions are satisfied by the Chen-Mangasarian smoothing and the Moreau envelope smoothing described in \cref{ex:chenmangasarian} and \cref{ex:moreau}. Under these circumstances, we obtain the uniform convergence \alert{\cref{assume:strongergradconsistency}}.%\jh{A concrete example is the hyperparameter learning model.. either put the example here or summarize everything in Discussions section later.} 
	%
	%
	%
	%
	%Often, the nonsmoothness of the lower-level function $g$ comes from the $y$ part. That is, $g(x,\cdot)$ may be nonsmooth for some $x$, but $g(\cdot,y)$ is a smooth function for all $y$. For instance, we may write $g(x,y) = g^{(1)}(x,y) + g^{(2)}(x,g^{(3)}(y)) + g^{(4)}(y)$ where $g^{(1)}$ and $g^{(2)}$ are smooth functions, but $g^{(3)}$ and $g^{(4)}$ are nonsmooth functions of $y$. In this case, we may consider smooth approximations for $g^{(3)}$ and $g^{(4)}$, say $g^{(3)}_{\mu}$ and $g^{(4)}_{\mu}$, to approximate $g$ as  $\gmu (x,y) = g^{(1)}(x,y) + g^{(2)}(x,\gmu ^{(3)}(y)) + \gmu ^{(4)}(y)$. Then $\nabla_x \gmu (x,y) = \nabla_x g^{(1)}(x,y) + \nabla_x g^{(2)}(x,\gmu ^{(3)}(y))$, and provided that $\gmu ^{(3)} $ converges uniformly to $g^{(3)}$ on a given domain, as required in the hypotheses of \cref{thm:smoothing_integral}, we get the uniform convergence assumption of \cref{thm:gradconsistency_integral}. In this example, note that $\nabla_x g(\cdot,\cdot)$ is continuous, and therefore we also achieve the upper semicontinuity assumption in \cref{thm:gradconsistency_integral}. 
\end{remark}

\begin{remark}
	\label{rem:uppersemicontinuity_partial_x_g}
	Upper semicontinuity of $\partial_x g$, \alert{as required in \cref{thm:gradconsistency_integral},} can be attained under any condition that will guarantee that if $(\xi,\xi')\in \partial g(\bar{x},\bar{y})$, then $\xi \in \partial_x g(\bar{x},\bar{y})$ for any $(\bar{x},\bar{y})$. For instance, the latter condition can be satisfied under Clarke regularity of $g$ (see \cite[Proposition 2.3.15]{Clarke83}) or convexity/concavity of $g$ with respect to $x$ (see \cite[Proposition 2.5.3]{Clarke83}). To see how the said condition implies upper semicontinuity of $\partial_x g$, note that given any $\varepsilon>0$, we have from upper semicontinuity of $\partial g$ that we can find $\delta>0$ such that $\partial g(x,y) \subseteq \partial g(\bar{x},\bar{y}) + \varepsilon B_1(0)$ for all $(x,y)\in B_{\delta}(\bar{x},\bar{y})$. From \cite[Proposition 2.3.16]{Clarke83}, we have $\partial_x g(x,y) \subset \pi_x \partial g(x,y)$. Finally, we get the desired result by noting that $\pi_x (\partial g(x,y)) \subset \pi_x (\partial g(\bar{x},\bar{y})+ \varepsilon B_1(0)) $.
\end{remark}

\begin{example}
	For the hyperparameter learning lower-level function $g$ given by \eqref{ex:hyperparameter_learning} where the $p_i$'s are not necessarily convex but may be nonsmooth, we may consider any smooth approximation $$\gmu (x,y) \coloneqq \ell_{\mu}(y) + \sum_{i=1}^n x_i p_{i,\mu}(y),$$
	where $p_{i,\mu} \to p_i$ uniformly on $Y\subseteq \Re^m$. When $p_i$ is convex, such uniform convergence required can be attained by using the Moreau envelope, or the Chen-Mangasarian smoothing. For nonconvex $p_i$ such as the $l_p$-norm with $p\in (0,1)$, uniform convergence holds for several examples as well (see \cite{ChenXiaojun12}). For this example of $g$ and $\gmu$, it is not difficult to verify that the assumptions required by \cref{thm:gradconsistency_integral} are all satisfied.  
\end{example}
With \cref{thm:gradconsistency_integral}, we immediately obtain gradient consistency as follows. 
\begin{theorem}
	\alert{Let $\bar{x}\in O_X$. Suppose that \cref{assume:quadred_4,assume:strongergradconsistency} hold}. Then the family of smooth approximations $\{ \alert{\vmuentrop} :\mu >0\}$ of $v$ satisfies the gradient consistent property at $\bar{x}$ under any of the following conditions:
	\begin{enumerate}[(a)]
		\item $\nabla_x g(x,y)$ exists and is continuous on $\Omega \times Y$; 
		\item $g(\cdot,y)$ is $\rho$-weakly concave function on $\Omega$ for any $y\in Y$; or
		\item \eqref{eq:partialformula} holds and $g$ is convex in $(x,y)$,
	\end{enumerate}
	where $\Omega$ is a neighborhood of $\bar{x}$.
	\label{cor:gradconsistency_integral}
\end{theorem}
\ifdefined\submit
\noindent {\textbf{Proof.}}
\else
\begin{proof}
\fi 
	Under any of the conditions in (a), (b), and (c), the set-valued mapping $\partial_x g$ is upper semicontinuous; see also \cref{rem:uppersemicontinuity_partial_x_g}. Hence, this theorem directly follows from \cref{thm:gradconsistency_integral,thm:danskin,thm:danskin_extended,thm:danskin_convex}.
\ifdefined\submit
\hfill \Halmos
\else 
\end{proof}
\fi

\subsection{Entropic regularization for unbounded constraint set}
\label{sec:entropy_noncompact}

We show that  when $Y$ is a closed unbounded set, the family of functions given by
\begin{equation}
\alert{\vmuentrop} (x) \coloneqq -\mu \ln \left( \int_{\Ymu} \exp \left( -\mu^{-1}\gmu(x,y)\right) dy \right),
\label{eq:smoothing_extended_linxuye_noncompact}
\end{equation}
is a smoothing function for the value function $v(x) = \max_{y\in Y} g(x,y)$. Throughout this section, we assume that $\{ \Ymu : \mu >0\}$ is a family of \textit{compact} subsets of $Y$ such that $\Ymu\nearrow Y$.

\begin{proposition}
	\label{prop:smoothing_entropy_noncompact}
	Suppose $\{\gmu: \mu >0\}$ is a family of smoothing functions of $g$ over $O_X\times O_Y$. For each $\mu>0$, $\alert{\vmuentrop}$ given by \eqref{eq:smoothing_extended_linxuye_noncompact} is continuously differentiable on $O_X$.
\end{proposition}
\ifdefined\submit
\noindent {\textbf{Proof.}}
\else
\begin{proof}
\fi 
	The claim can be proved by replacing $Y$ in the proof of \cref{prop:smoothing_entropy} with $\Ymu$. We note that the gradient of $\alert{\vmuentrop}$ is given by 
	\begin{equation*}
	\nabla \alert{\vmuentrop} (x) = \int_{\Ymu} \alpha_{\mu}(x,y) \nabla_x \gmu (x,y) dy ,
	\label{eq:vmu_integral_noncompact}
	\end{equation*}
	where 
	\begin{equation}
	\alpha_{\mu} (x,y) \coloneqq \frac{\exp(-\mu^{-1}\gmu (x,y)) }{\int_{\Ymu} \exp(-\mu^{-1}\gmu (x,z))  dz}.
	\label{eq:alpha_mu_noncompact}
	\end{equation}
\ifdefined\submit
\hfill \Halmos
\else 
\end{proof}
\fi 

\alert{Analogous to \cref{prop:ennamcq_entrop_compact}, we have the following result. As will be shown in \cref{thm:smoothing_integral_noncompact}, the condition $ \mu \ln \vol(Y_{\mu}) \to 0 $ as $\mu\to 0$ (that is, \cref{assume:controlledvolume}) is necessary for obtaining a smooth approximation of $ v $. Consequently, the condition $ \mu \ln \vol(Y_{\mu}) < \epsilon $ required in the proposition below will be satisfied under \cref{assume:controlledvolume} for sufficiently small $\mu$.}
\alert{\begin{proposition}\label{prop:ennamcq_entrop_noncompact}
		Suppose that \cref{assume:quadreg}(b) holds, and let $v_{\mu}=\vmuentrop$ in \eqref{eq:bilevel_nlp_formulation_epsilon_mu} with $\epsilon>0$. For any $\mu>0$ such that $\mu \ln \vol (Y_{\mu}) < \epsilon$, ENNAMCQ holds at any point $(\bar{x},\bar{y})\in X\times Y$. 
		%		, that is, either $g_{\mu}(\bar{x},\bar{y})-\vmuentrop (\bar{x})<\epsilon$ or $g_{\mu}(\bar{x},\bar{y})-\vmuentrop (\bar{x})\geq \epsilon$ but $0\notin \nabla g_{\mu}(\bar{x},\bar{y}) - \nabla \vmuentrop(\bar{x})\times \{0\} + N_{X\times Y}(\bar{x},\bar{y}). $
\end{proposition}}

\alert{To prove an analogue of \cref{lemma:integral_lower_bound}, we require the uniform level-boundedness assumption in \cref{assume:quadred_3}. }
\begin{lemma}
	\label{lemma:integral_lower_bound_noncompact}
	\alert{Suppose that \cref{assume:quadred_3} holds}. Given any $\bar{x}\in O_X$, there exist a neighborhood $\Omega$ of $\bar{x}$ and a positive number $\mu_0>0$ such that 
	\begin{equation}
	\max_{y\in Y_{\mu_0}} \exp (-g(x,y)) = \max_{y\in Y} \exp (-g(x,y)) \qquad \forall x\in \cl (\Omega).
	\label{eq:max_are_equal_noncompact}
	\end{equation}
	Moreover, for any $\tau \in (0,1)$, there exists $\delta>0$ such that 
	\begin{equation}
	\tau (\mu \vol(Y_{\mu_0}))^{\mu} \max _{y\in Y} \exp (-g(x,y)) \leq \left( \int_{\Ymu} \exp \left( -\mu^{-1}g(x,y)\right)dy \right) ^{\mu} \quad \forall x\in \cl (\Omega), ~\forall \mu\in (0,\delta). 
	\label{eq:integral_lower_bound_noncompact}
	\end{equation}
\end{lemma}

\ifdefined\submit
\noindent {\textbf{Proof.}}
\else
\begin{proof}
\fi 
	From the proof of \cref{thm:danskin_extended_noncompact}, we know that $\bigcup_{x\in  \cl (B_{\alpha}(\bar{x}))} S(x)$ is a bounded set for some $\alpha<1$. Therefore, we may choose $\mu_0$ small enough so that $Y_{\mu_0}$ contains $\bigcup_{x\in \cl (\Omega)} S(x)$ with $\Omega =  B_{\alpha}(\bar{x})$. Hence, \eqref{eq:max_are_equal_noncompact} holds. Now, let $\tau\in (0,1)$. We have from the proof of \cref{lemma:integral_lower_bound} that there exists a sufficiently large $N_0$ such that  	
	\begin{equation}
	\ds \max_{j=1,2,\dots,k} ~~\min_{y\in (Y_{\mu_0})_j} \exp (-g(x,y)) > \tau \cdot \max_{y\in Y_{\mu_0}} \exp (-g(x,y)) \quad \forall x\in \cl (\Omega)
	\label{eq:maxminmaxinequality_noncompact}
	\end{equation}
	holds for all integers $k\geq N_0$ and any partition $\{(Y_{\mu_0})_1,(Y_{\mu_0})_2\dots,(Y_{\mu_0})_k\}$ of $Y_{\mu_0}$ with $\ds \max_{j=1,2,\dots,k} \diam (Y_{\mu_0})_j\to 0$ as $k\to \infty$. Choose $\delta < 1/N_0$ and \alert{a partition $\{(Y_{\mu_0})_1,(Y_{\mu_0})_2\dots,(Y_{\mu_0})_{N_0}\}$ of $Y_{\mu_0}$ such that $\vol ((Y_{\mu_0})_i) = \vol ((Y_{\mu_0})_j)$ for any $i,j$}. Then, for any $\mu \in (0,\delta)$,
	\begin{align*}
	\left( \int_{\Ymu} \exp \left( -\mu^{-1}g(x,y)\right)dy \right) ^{\mu} & \geq \left( \int_{Y_{\mu_0}} \exp \left( -\mu^{-1}g(x,y)\right)dy \right) ^{\mu} \\
	& \geq \left( \frac{\vol(Y_{\mu_0})}{N_0} \sum_{j=1}^{N_0} \min_{y\in (Y_{\mu_0})_j} \exp \left( -\mu^{-1}g(x,y)\right) \right)^{\mu} \\ 
	& = \left( \frac{\vol(Y_{\mu_0})}{N_0} \sum_{j=1}^{N_0} \left( \min_{y\in (Y_{\mu_0})_j} \exp \left( -g(x,y)\right) \right)^{\mu^{-1}} \right)^{\mu} \\ 
	& \geq  \left( \frac{\vol(Y_{\mu_0})}{N_0}   \left( \max _{j=1,\dots,N_0}\min_{y\in (Y_{\mu_0})_j} \exp \left( -g(x,y)\right) \right) ^{\mu^{-1}}\right)^{\mu} \\
	& \geq \tau \left( \frac{\vol(Y_{\mu_0})}{N_0}\right)^{\mu}  \max_{y\in Y_{\mu_0} }\exp (-g(x,y)) \\
	& \geq \tau  \vol(Y_{\mu_0})^{\mu} \mu^{\mu}\max_{y\in Y_{\mu_0} }\exp (-g(x,y))\\
	& =\tau  \vol(Y_{\mu_0})^{\mu} \mu^{\mu}\max_{y\in Y }\exp (-g(x,y)),
	\end{align*}
	where the fourth inequality holds by \eqref{eq:maxminmaxinequality_noncompact}, and the last equality holds by \eqref{eq:max_are_equal_noncompact}. This completes the proof. 
\ifdefined\submit
\hfill \Halmos
\else 
\end{proof}
\fi 

We now show that the family $\{\vmuentrop : \mu > 0\}$ constitutes a smooth approximation in the sense of \cref{defn:smoothing}. \alert{To this end, we impose \cref{assume:controlledvolume} on the sets $\{ \Ymu : \mu > 0 \}$. Note that a simple choice that satisfies this condition is $\Ymu \coloneqq Y \cap B_{r_{\mu}}(0)$, where $r_{\mu} \nearrow +\infty$ as $\mu \searrow 0$ and $\vol(B_{r_{\mu}}(0)) \leq \mu^{-p}$ for any $p > 0$. 
\begin{assumption}
	\label{assume:controlledvolume}
	The family $\{\Ymu: \mu>0\}$ is chosen so that $\mu \ln \vol(\Ymu)\to 0$ as $\mu\to 0$
\end{assumption}}
 
\begin{theorem}
	\label{thm:smoothing_integral_noncompact}
	\alert{Suppose that \cref{assume:quadred_3,assume:quadred_4,assume:controlledvolume} hold}. Then $\{\alert{\vmuentrop} : \mu >0\}$ is a family of smooth approximations for the value function over $O_X$
\end{theorem}

\ifdefined\submit
\noindent {\textbf{Proof.}}
\else
\begin{proof}
\fi 
	We only need to show that for any $\bar{x}\in O_X$, $\alert{\vmuentrop}(x)\to v(\bar{x})$ as $\mu\to 0$ and $x\to \bar{x}$. To this end, let $\varepsilon >0$ be given, and let $\tau \in (\exp (-\varepsilon),1)$. By \cref{lemma:integral_lower_bound_noncompact} there exist a neighborhood $\Omega$ of $\bar{x}$ and positive numbers $\mu_0$ and $\delta<\mu_0$ such that \eqref{eq:max_are_equal_noncompact} and \eqref{eq:integral_lower_bound_noncompact} hold for all $\mu\in (0,\delta)$. Meanwhile, arguing as in \cite[Theorem 1]{FangWu96}, we have
	\[ \left( \int_{\Ymu} \exp \left( -\mu^{-1}g(x,y)\right) dy\right) ^{\mu}\leq \vol(\Ymu)^{\mu} \max_{y\in \Ymu} \exp (-g(x,y)), \quad \forall x\in O_X.\]
	Together with \eqref{eq:integral_lower_bound_noncompact}  and noting that $Y_{\mu} \supseteq Y_{\mu_0}$, we have 
	\[\tau (\mu \vol(Y_{\mu_0}))^{\mu} \max _{y\in Y} \exp (-g(x,y)) \leq \left( \int_{\Ymu} \exp \left( -\mu^{-1}g(x,y)\right) dy\right) ^{\mu}\leq \vol(\Ymu)^{\mu} \max_{y\in Y} \exp (-g(x,y)) \]
	for all $x\in \cl (\Omega)$. Since the logarithmic function is monotonic, then 
	\begin{multline*}
	- \mu \ln \vol(\Ymu) -  \max_{y\in Y} (-g(x,y)) \leq -\mu \ln \left( \int_{\Ymu} \exp \left( -\mu^{-1}g(x,y)\right) dy\right) \\ \leq -\ln \tau - \mu \ln  (\mu \vol(Y_{\mu_0}))-  \max _{y\in Y} (-g(x,y)).
	\end{multline*}
	Then
	\begin{equation*}
	v(x) - \mu \ln \vol(\Ymu) \leq -\mu \ln \left( \int_{\Ymu} \exp \left( -\mu^{-1}g(x,y)\right) dy\right)  \leq v(x) - \mu \ln  (\mu \vol(Y_{\mu_0})) + \varepsilon.
	\end{equation*}
	\alert{By \cref{assume:controlledvolume}, we} may choose $\delta>0$ to be sufficiently small enough so that for all $\mu\in (0,\delta)$, we have $ \mu \ln \vol(\Ymu)<2\varepsilon$ and $ - \mu \ln  (\mu \vol(Y_{\mu_0}))< \varepsilon $. Hence, 
	\begin{equation*}
	\left| -\mu \ln \left( \int_{\Ymu}  \exp \left( -\mu^{-1}g(x,y)\right) dy\right) - v(x) \right| < 2\varepsilon,
	\label{eq:bounds_v(x)_noncompact}
	\end{equation*}
	valid for all $x\in \cl(\Omega)$ and $\mu\in (0,\delta)$. The rest of the proof follows from the same arguments as in \cref{thm:smoothing_integral}. 
	
	%	On the other hand, invoking the uniform convergence of $\gmu$ to $g$ on $\cl (\Omega)\times Y$, we may choose $\delta$ to be small enough so that \eqref{eq:bounds_v(x)_noncompact} holds and $|\gmu (x,y) - g(x,y)|<\varepsilon$ for all $\mu \in (0,\delta)$ and $(x,y)\in \cl (\Omega)\times Y$. Consequently, we have 
	%	\begin{equation}
	%	\left| -\mu \ln \left( \int_Y \exp \left( -\mu^{-1}\gmu (x,y)\right) dy\right) + \mu \ln \left( \int_Y \exp \left( -\mu^{-1}g(x,y)\right) dy\right) \right| <  \varepsilon
	%	\end{equation}
	%	for all $x\in \cl (\Omega)$ and $\mu \in (0,\delta)$. Together with \eqref{eq:bounds_v(x)_noncompact} and the definition of $\alert{\vmuentrop}$, we have 
	%	\begin{equation}
	%	| \alert{\vmuentrop} (x) - v(x) | < 3\varepsilon \qquad \forall x\in \cl (\Omega), ~\forall \mu \in (0,\delta) .
	%	\label{eq:bounds_vmu_noncompact}
	%	\end{equation}
	%	By the continuity of the value function (see \cref{thm:uniformlevelbounded}(a)) at $\bar{x}\in \cl (\Omega)$, we may choose $\delta>0$ so that $|v(x) - v(\bar{x})|< \varepsilon$ for all $x\in \cl (\Omega)$ with $\|x - \bar{x}\|<\delta$. Together with \eqref{eq:bounds_vmu_noncompact}, we have $|\alert{\vmuentrop} (x) - v(\bar{x}) |< 4\varepsilon$ for all $\mu\in (0,\delta)$ and $x\in \cl (\Omega)$ with $\|x-\bar{x}\|<\delta$. This completes the proof. 
\ifdefined\submit
\hfill \Halmos
\else 
\end{proof}
\fi

As shown above, managing the unboundedness of the constraint set $Y$ to obtain a smoothing function of the value function can be achieved by assuming uniform level-boundedness of $g$ and restricting the region of integration to $\Ymu$. In particular, we need the $\Ymu$'s to be compact sets whose volumes increase in a controlled magnitude so that $\mu \ln (\vol (\Ymu))\to 0$. However, these are not sufficient to prove gradient consistency of the smoothing family, and the main difficulty lies on the unbounded nature of $Y$. Integrable dominating functions would be required to handle limits of sequences of integrals as shown in the following results. We first note the following counterpart of \cref{lemma:alpha_outside_solutionset}.

\begin{lemma}
	 \alert{Suppose that \cref{assume:quadred_3,assume:quadred_4,assume:controlledvolume} hold. If $\bar{x}\in O_X$ and $y'\in Y\setminus S(\bar{x})$, then} $\alpha_{\mu}(x,y')\to 0$ as $(x,\mu)\to (\bar{x},0)$, where $\alpha_{\mu}$ is given by \eqref{eq:alpha_mu_noncompact}. 
	\label{lemma:alpha_outside_solutionset_noncompact}
\end{lemma}

\ifdefined\submit
\noindent {\textbf{Proof.}}
\else
\begin{proof}
\fi 
	The proof is exactly identical to that of \cref{lemma:alpha_outside_solutionset},
\ifdefined\submit
\hfill \Halmos
\else 
\end{proof}
\fi 

To prove a counterpart of \cref{lemma:alpha_outside_solutionset_Integral} for the unbounded case, we need to assume the existence of an integrable dominating function $h$ for $\alpha_{\mu}(x,\cdot)$ over $Z$ for all $(x,\mu)$ in a neighborhood of $(\bar{x},0)$. Note that in the case of \cref{lemma:alpha_outside_solutionset_Integral}, we have shown in the proof that uniform convergence of $\{\gmu\}$ and the compactness of $Y$ lead us to obtain the constant function $h\equiv 1$ as the desired dominating function, which is indeed integrable over the compact set $Z$. 
\alert{\begin{assumption}
	\label{assume:dominatingfunction}
	At $\bar{x}\in O_X$, we have that for any closed subset $Z$ of $Y\setminus S(\bar{x})$,  there exists an integrable function $h:Z\to \Re$ such that $\alpha_{\mu}(x,y) \leq h(y)$ for all $y\in Z$ and for all $(x,\mu)$ in some neighborhood of $(\bar{x},0)$.
\end{assumption}}
\begin{lemma}
	\label{lemma:alpha_outside_solutionset_Integral_noncompact}
	\alert{Let $\bar{x}\in O_X$. If \cref{assume:quadred_3,assume:quadred_4,assume:controlledvolume,assume:dominatingfunction} hold, then} $\int_{Z} \alpha_{\mu}(x,y) dy\to 0$ as $(x,\mu) \to (\bar{x},0)$. 
\end{lemma}

\ifdefined\submit
\noindent {\textbf{Proof.}}
\else
\begin{proof}
\fi 
	This is a direct consequence of the Lebesgue dominated convergence theorem and  \cref{lemma:alpha_outside_solutionset_noncompact}.
\ifdefined\submit
\hfill \Halmos
\else 
\end{proof}
\fi 

The following is another technical lemma needed to guarantee that $\{ \nabla \alert{\vmuentrop} (x)\}$ is bounded over all $(x,\mu)$ in a neighborhood of $(\bar{x},0)$. 
\alert{\begin{assumption}
	\label{assume:dominatingfunction2}
	At $\bar{x}\in O_X$, we have that for any closed subset $Z$ of $Y\setminus S(\bar{x})$, there exists an integrable function $H:Z\to \Re$ such that $\alpha_{\mu}(x,y)\norm{\nabla_x \gmu (x,y)} \leq H(y)$ for all $y\in Z$ and for all $(x,\mu)$ in some neighborhood of $(\bar{x},0)$.
\end{assumption}}
\begin{lemma}
	\alert{Let $\bar{x}\in O_X$. Suppose that \cref{assume:quadred_3,assume:quadred_4,assume:strongergradconsistency,assume:controlledvolume,assume:dominatingfunction,assume:dominatingfunction2} hold,  and that there exists a neighborhood $\Omega$ of $\bar{x}$ such that $\partial_x g(\cdot,\cdot)$ is upper semicontinuous on $\Omega\times O_Y$}. Then there exists $\delta'>0$ such that $\{ \nabla \alert{\vmuentrop} (x) : \mu\in (0,\delta') , x\in \cl (B_{\delta'}(\bar{x}))\}$ is bounded. 
	%That is, there exists $M>0$ such that 
	%\begin{equation}
	%\norm{\alert{\vmuentrop}(x)} \leq M \quad \forall \mu\in (0,\delta'), ~\forall x\in \cl (B_{\delta'}(\bar{x})).
	%\label{eq:bounded_gradient_vmu}
	%\end{equation}
	\label{lemma:bounded_gradient_vmu}
\end{lemma}

\ifdefined\submit
\noindent {\textbf{Proof.}}
\else
\begin{proof}
\fi 
	Let $\Theta$ be a bounded open set that properly contains the compact set $S(\bar{x})$, and let $Z\coloneqq Y\setminus \Theta$, which is a closed set contained in $Y\setminus S(\bar{x})$. Meanwhile, we know from the uniform convergence of $\dist (\nabla_x \gmu (x,\cdot), \partial_x g(x,\cdot)) $ on $Y$ that there exists $\delta>0$ such that 
	\begin{equation}
	\dist ( \nabla_x \gmu (x,y), \partial_x g(x,y)) < 1 \quad \forall \mu \in (0,\delta) , ~~\forall (x,y)\in \cl (B_{\delta}(\bar{x}))\times Y.
	\label{eq:uniconv_noncompact}
	\end{equation}
	We may argue as in the first part of the proof of \cref{thm:gradconsistency_integral} to conclude that the set 
	\begin{equation}
	\left\lbrace \int_{\Theta \cap Y} \alpha_{\mu}(x,y) \| \nabla_x \gmu (x,y)\| dy ~:~   x\in \cl (B_{\delta}(\bar{x})), ~\mu\in (0,\delta) \right\rbrace 
	\label{eq:integraloverYmu0}
	\end{equation}
	is bounded, noting the set $\cl (\Theta)\cap Y$ is compact.  From the upper semicontinuity of $\partial_x g(\cdot,\cdot)$, we see that $\bigcup _{x\in \cl (B_{\delta}(\bar{x}))} \partial_x g(x,y)$ is bounded for any given $y\in Z$. Together with \eqref{eq:uniconv_noncompact} and triangle inequality, we know that $\{ \norm{\nabla_x g(x,y)} : x\in \cl (B_{\delta}(\bar{x})), ~\mu\in (0,\delta)\}$ is bounded for a fixed $y$. Hence, using \cref{lemma:alpha_outside_solutionset_Integral_noncompact}, we conclude that $\alpha_{\mu}(x,y) \norm{\nabla_x g(x,y)} $ converges pointwise to $0$ as $(x,\mu)\to(\bar{x},0)$ for each $y\in Z$. Since an integrable dominating function exists by our hypothesis, we conclude by the Lebesgue dominated convergence theorem that 
	\begin{equation}
	\int_{Z} \alpha_{\mu}(x,y) \|\nabla_x \gmu (x,y)\|  dy \to 0 \quad \text{as} \quad (x,\mu)\to (\bar{x},0).
	\label{eq:integralgoestozero_noncompact}
	\end{equation}
	Finally, since 
	\begin{align*}
	\norm{\nabla \alert{\vmuentrop} (x)} & \leq \int_{\Ymu} \alpha_{\mu}(x,y) \norm{\nabla_x \gmu (x,y)} dy\\ 
	& = \int_{\Theta\cap Y} \alpha_{\mu}(x,y) \norm{\nabla_x \gmu (x,y)} dy  + \int_{\Ymu\setminus (\Theta \cap Y)} \alpha_{\mu}(x,y) 
	\norm{\nabla_x \gmu (x,y)} dy  \\
	& \leq \int_{\Theta\cap Y} \alpha_{\mu}(x,y) \norm{\nabla_x \gmu (x,y)} dy  + \int_{Z} \alpha_{\mu}(x,y) \norm{\nabla_x \gmu (x,y)} dy ,
	\end{align*}
	the boundedness of \eqref{eq:integraloverYmu0} and the limit \eqref{eq:integralgoestozero_noncompact} imply the boundedness of $ \norm{ \nabla \alert{\vmuentrop} (x)}$ on a neighborhood of $(\bar{x},0)$, thus completing the proof.
\ifdefined\submit
\hfill \Halmos
\else 
\end{proof}
\fi 

\begin{lemma} 	\alert{Let $\bar{x}\in O_X$. Suppose that \cref{assume:quadred_3,assume:quadred_4,assume:strongergradconsistency,assume:controlledvolume,assume:dominatingfunction,assume:dominatingfunction2} hold,  and that there exists a neighborhood $\Omega$ of $\bar{x}$ such that $\partial_x g(\cdot,\cdot)$ is upper semicontinuous on $\Omega\times O_Y$.} Then 
	\[\emptyset \neq \limsup_{x\to \bar{x}, \mu\searrow 0}\nabla \alert{\vmuentrop} (x) \subseteq P(\bar{x}), \]
	where $P(\bar{x})$ is given in \alert{\eqref{eq:P(x)}}. 
	\label{thm:gradconsistency_integral_noncompact}
\end{lemma}
\ifdefined\submit
\noindent {\textbf{Proof.}}
\else
\begin{proof}
\fi 
	That $\limsup_{x\to \bar{x}, \mu\searrow 0}\nabla \alert{\vmuentrop} (x)$ is nonempty follows from \cref{lemma:bounded_gradient_vmu}. To prove the inclusion, suppose that $\nabla v_{\mu_k}(x^k)\to \xi\in\Re^n$ for some sequences $x^k\to \bar{x}$ and $\mu_k\searrow 0$ as $k\to \infty$. We need to show that $\xi \in P(\bar{x})$. As in \cref{thm:gradconsistency_integral}, it suffices to show that for any given $\varepsilon>0$, \alert{there exists an element of $\xi^k \in P(\bar{x})$ such that $\norm{\nabla v_{\mu_{k}}(x^{k})-\xi^k}<\varepsilon$,} if $k$ is chosen sufficiently large. 
%	That is, there exists $\xi^k\in P(\bar{x})$ such that $\| \nabla v_{\mu_{k}}(x^{k}) - \xi^k\| <\varepsilon$. 
	Following similar arguments in \cref{thm:gradconsistency_integral},  we know that we can find points $y_1,y_2,\dots, y_r\in S(\bar{x})$ and numbers $\delta_1,\delta_2,\dots,\delta_r>0$ such that $S(\bar{x})\subseteq \bigcup_{i=1}^r B_{\delta_i}(y_i)$ and
	\begin{equation*}
	\partial_x g(x',y') \subseteq \partial_x g(\bar{x},y_i) + \frac{\varepsilon}{3} B_1(0)\quad \forall x'\in B_{\delta_i}(\bar{x}), ~y'\in B_{\delta_i}(y_i),
	\label{eq:usc_noncompact}
	\end{equation*}
	for $i=1,\dots, r$. Denote $Y_1 \coloneqq Y\cap B_{\delta_1}(y_1)$, and $Y_i  \coloneqq  Y \cap \left(B_{\delta_i} (y_i) \setminus \bigcup_{j=1}^{i-1}B_{\delta_j}(y_j) \right)$ for $i=2,\dots, r$. Moreover, we let $y_{r+1}$ be an arbitrary element of $S(\bar{x})$, and we denote $Y_{r+1}^k \coloneqq Y_{\mu_k}\setminus \bigcup_{i=1}^r B_{\delta_i}(y_i)$.
	
	From the uniform convergence of $\dist (\nabla_x \gmu (x,\cdot), \partial_x g(x,\cdot)) $ on $Y$, we can find a $\delta>0$ such that $\delta<\min \{ \delta_1,\dots, \delta_r\}$ and 
	\begin{equation}
	\dist ( \nabla_x \gmu (x,y), \partial_x g(x,y)) < \frac{\varepsilon}{3} \quad \forall \mu \in (0,\delta)  ~\text{and}~\forall (x,y)\in \cl (B_{\delta}(\bar{x}))\times Y.
	\label{eq:uniconvergence_derivative_noncompact}
	\end{equation}
	Setting $\Theta \coloneqq \bigcup_{i=1}^r B_{\delta_i}(y_i)$ and $Z \coloneqq Y\setminus \Theta$, we know that \eqref{eq:integralgoestozero_noncompact} holds since $\Theta$ contains $S(\bar{x})$. Hence, we may choose a smaller $\delta$ if necessary so that  \eqref{eq:uniconvergence_derivative_noncompact}  holds and at the same time,
	\begin{equation}
	\int_{Z}\alpha_{\mu} (x,y) \norm{\nabla_x g(x,y)}dy < \frac{\varepsilon}{6} \qquad \forall \mu\in (0,\delta),~ \forall  x\in \cl(B_{\delta}(\bar{x})).
	\label{eq:integral_epsilon_6}
	\end{equation}
	Choosing a smaller $\delta$ if needed, we have from \cref{lemma:alpha_outside_solutionset_Integral_noncompact} that 
	\begin{equation}
	\int_Z \alpha_{\mu}(x,y) < \frac{\varepsilon}{6M} \qquad \forall \mu\in (0,\delta), ~\forall x\in \cl(B_{\delta}(\bar{x})),
	\label{eq:integral_epsilon_6M}
	\end{equation}
	where $M \coloneqq \max \{\|\xi\|: \xi \in \partial_x g(\bar{x},y_{r+1}) \}$, which is finite due to the boundedness of  $\partial_x g(\bar{x},y_{r+1}) $. 
	
	Now, let us choose $k$ to be large enough so that $x^k\in B_{\delta}(\bar{x})$ and $\mu_k \in (0,\delta)$. Define $\xi^k$ as 
	\begin{equation*}
	\xi^k \coloneqq \sum_{i=1}^{r}\int_{Y_i} \alpha_{\mu_k}(x^k,y) \Pi_{\partial_xg(\bar{x},y_i)} \left( \varphi^{(i)}_k(y)\right)dy + \int_{Y_{r+1}^k} \alpha_{\mu_k}(x^k,y) \Pi_{\partial_xg(\bar{x},y_{r+1})} \left( \varphi^{(r+1)}_k(y)\right)dy ,
	\end{equation*}
	where $\varphi^{(i)}_k$ ($i=1,\dots,r$) is \alert{given by \eqref{eq:varphi_i}}, and $\varphi^{(r+1)}_k$ is a single-valued mapping from $Y_{r+1}^k$ to $\partial_x g(x^k,y)$. Note that since $Y_{r+1}^k \cup \left(\bigcup_{i=1}^{r} Y_i \right) = Y_{\mu_k}$, $\int_{Y_{\mu_k}} \alpha_{\mu_k} (x^k,y)dy = 1$, and the range of $\Pi_{\partial_xg(\bar{x},y_i)}$ is contained in $P(\bar{x})$, it follows that $\xi^k$ is indeed in $P(\bar{x})$ by Jensen's inequality. 
	
	We show that $\| \nabla v_{\mu_{k}}(x^{k}) - \xi^k\| <\varepsilon$. 
%	As before, we denote by $E_i$ the set $E_i \coloneqq \{y\in Y_i : g~\text{is not differentiable at }(x^k,y) \} $ for each $i=1,\dots, r$, and we note that each $E_i$ has zero Lebesgue measure. 
	\alert{We can follow the same arguments in \cref{thm:gradconsistency_integral} to show that 
	\begin{equation}
	\norm{ \nabla_x g_{\mu_k}(x^k,y) - \Pi_{\partial_xg(\bar{x},y_i)} \left( \varphi^{(i)}_k(y)\right)} < \frac{2\varepsilon}{3}, \quad \forall y\in Y_i, ~i\in \{1,\dots,r\}.\label{eq:dist_S(x)_integral_noncompact}
	\end{equation}}
	On the other hand, we also have
	\begin{align}
	&\norm{\int_{Y_{r+1}^k} \alpha_{\mu_k}(x^k,y) \left(\nabla_x g_{\mu_k}(x^k,y)- \Pi_{\partial_xg(\bar{x},y_{r+1})} \left( \varphi^{(r+1)}_k(y)\right)\right)dy } \notag \\
	&\quad \leq  \int_{Y_{r+1}^k} \alpha_{\mu_k}(x^k,y)  \norm{\nabla_x g_{\mu_k}(x^k,y)}dy + \int_{Y_{r+1}^k} \alpha_{\mu_k}(x^k,y)  \norm{\Pi_{\partial_xg(\bar{x},y_{r+1})} \left( \varphi^{(r+1)}_k(y)\right)} dy \notag \\
	&\quad \leq  \int_{Z} \alpha_{\mu_k}(x^k,y)  \norm{\nabla_x g_{\mu_k}(x^k,y)}dy + \int_{Z} \alpha_{\mu_k}(x^k,y)  \norm{\Pi_{\partial_xg(\bar{x},y_{r+1})} \left( \varphi^{(r+1)}_k(y)\right)} dy \notag \\
	& \quad \overset{\eqref{eq:integral_epsilon_6}}{<} \frac{\varepsilon}{6} + M \int_{Z} \alpha_{\mu_k}(x^k,y)  dy \notag \\  
	& \quad \overset{\eqref{eq:integral_epsilon_6M}}{<}\frac{\varepsilon}{3}, \label{eq:outside_S(x)_integral_noncompact}
	\end{align}
	where second inequality holds since $Y_{r+1}^k \subseteq Z$ for any $k$. Finally, we can use \eqref{eq:outside_S(x)_integral_noncompact} and \eqref{eq:dist_S(x)_integral_noncompact} to show that $\| \nabla v_{\mu_k}(x^k) - \xi^k\|<\varepsilon$ by arguing as in \eqref{eq:xi_epsilon_away}. 
\ifdefined\submit
\hfill \Halmos
\else 
\end{proof}
\fi

We finally show that $\vmu$ satisfies the gradient consistent property.

\begin{theorem}
	\alert{Let $\bar{x}\in O_X$. Suppose that \cref{assume:quadred_3,assume:quadred_4,assume:strongergradconsistency,assume:controlledvolume,assume:dominatingfunction,assume:dominatingfunction2} hold}. Then the family of smooth approximations $\{ \vmu :\mu >0\}$ of $v$ satisfies the gradient consistent property at $\bar{x}$ under any of the following conditions:
	\begin{enumerate}[(a)]
		\item $\nabla_x g(x,y)$ exists and is continuous on $\Omega\times O_Y$; or
		\item $g(\cdot,y)$ is $\rho$-weakly concave function on $\Omega$ for any $y\in Y$; or
		\item \eqref{eq:partialformula} holds and $g$ is convex in $(x,y)$, 
	\end{enumerate}
	where $\Omega$ is a neighborhood of $\bar{x}$
	\label{cor:gradconsistency_integral_noncompact}
\end{theorem}
\ifdefined\submit
\noindent {\textbf{Proof.}}
\else
\begin{proof}
\fi 
	As noted in the proof of \cref{cor:gradconsistency_integral}, $\partial_x g$ is upper semicontinuous under any of the conditions (a), (b) and (c). The result is then a direct consequence of \cref{thm:gradconsistency_integral_noncompact,thm:danskin_noncompact,thm:danskin_extended_noncompact,thm:danskin_convex}.
\ifdefined\submit
\hfill \Halmos
\else 
\end{proof}
\fi 

\section{Concluding Remarks \alert{and Future Directions}}\label{sec:conclusion}

We have presented two techniques for deriving smooth approximations of the value function corresponding to a nonsmooth function. One approach is based on quadratic regularization, which is applicable when the lower-level objective is convex in the decision variable of the lower-level problem. In the absence of this convexity structure, we proposed an alternative method via entropic regularization, which only requires Lipschitz continuity. While smoothing approaches are often employed for the practical purpose that optimization toolboxes usually require smoothness, it is equally important to theoretically investigate whether or not the smoothing algorithm will converge to a stationary point. To this end, we also proved the gradient consistent property for the proposed smoothing techniques, which is a central requirement for the convergence of smoothing algorithms. This paper provides the theoretical framework for smoothing-based methods for bilevel optimization problems. 

\alert{In future work, we aim to explore several directions. First, we will further investigate the possibility of relaxing the Lipschitz continuity assumption on the lower-level objective function. The existence of dominating functions, which was crucial for entropic regularization under unbounded constraint sets, also deserve deeper investigation. Another important direction is to extend our smoothing frameworks to settings where the lower-level constraint set is not fixed but instead depends on the upper-level variable \( x \). A preliminary approach based on quadratic regularization combined with penalization has been proposed in \citep{Liu23}; however, it remains an open question whether the resulting approximation qualifies as a smoothing family in the sense of \cref{defn:smoothing} and whether it satisfies gradient consistency. In connection with this, we also note the recent work \citep{xu2024enhanced}, which constructs smoothing functions for the KKT solution mapping of the lower-level problem with inequality constraints (assumed differentiable with respect to the lower-level variable) and analyzes gradient consistency. While related, our focus is distinct $-$ we aim to develop smoothing functions for the \textit{value function} of the lower-level problem, rather than its KKT solution mapping. Finally, a key part of our near-future work is the development of efficient smoothing algorithms grounded in the theoretical framework established in this paper. In this regard, we note that several recent works \citep{LamparielloSagratella20,Liu23,lu2023slm,pedregosa16,YeYuanZengZhang23} have introduced algorithms that rely on gradient approximations of the lower-level constraint functions, rather than exact evaluations. This line of work offers a promising direction to enhance computational efficiency. We plan to explore how such approximation-based strategies can be effectively combined with our framework to design practical algorithms with provable convergence guarantees.}

\ifdefined\submit
\ACKNOWLEDGMENT{
The authors are grateful to the associate editor and the two anonymous referees for their thoughtful and constructive comments on an earlier version of this paper. Their suggestions greatly improved the presentation of the results and the overall quality of the manuscript.
}	
\else 
\fi 

\ifdefined\submit 
\bibliographystyle{informs2014} % outcomment this and next line in Case 1
\else 
\fi 

\bibliography{bibfile} % if more than one, comma separated

	\end{document}